\documentclass[a4paper,preprint,10pt,maxbibnames=20,giveninits=true,sort&compress]{elsarticle}

\usepackage[bottom=2cm,left=2cm,right=2cm,top=2cm]{geometry}
\usepackage{amssymb}
\usepackage{xcolor}
\usepackage{float}
\usepackage[utf8]{inputenc}
\usepackage[T1]{fontenc}
\usepackage{amssymb,amsmath,amsthm}
\usepackage{graphicx}
\usepackage[colorlinks=false, pdfborder={0 0 0}]{hyperref}
\usepackage{csquotes}
\usepackage{subcaption}
\usepackage{accents}
\usepackage{textcomp}
\usepackage{import}
\usepackage{caption,graphicx,newfloat}
\usepackage[symbol]{footmisc}

\DeclareFloatingEnvironment[
fileext=lob,
listname={List of Boxes},
name=Box,
placement=htp,
]{BOX}

\newcommand{\dfracderiv}[2]{\dfrac{\mathrm{d}\,#1}{\mathrm{d}\,#2}}
\newcommand{\norm}[1]{\left\Vert#1\right\Vert}
\newcommand{\dfracderivpartial}[2]{\dfrac{\partial\,#1}{\partial\,#2}}
\newcommand{\matr}[2][{}]{\underline{#2}^{#1}}
\newcommand{\Matr}[2][{}]{\underline{\underline{#2}}^{#1}}
\newcommand{\vect}[2][{}]{\underline{#2}^{#1}}
\newcommand{\tens}[2][{}]{\undertilde{#2}^{#1}}

\newcommand{\indexnot}[3][{}]{{#2}^{#1}_{#3}}
\newcommand{\transp}[1]{{#1}^\mathrm{T}}
\newcommand{\invers}[1]{{#1}^{-1}}
\newcommand{\transpinv}[1]{{#1}^{-\mathrm{T}}}
\newcommand{\nodal}[1]{\hat{#1}}
\newcommand{\amplitudal}[1]{\Tilde{#1}}
\newcommand{\MStretch}{U}
\newcommand{\MDefgrad}{F}
\newcommand{\MJacobian}{J}
\newcommand{\MRateofdef}{D}

\newcommand{\MSpin}{W}
\newcommand{\MRotation}{R}
\newcommand{\MStress}{\Sigma}
\newcommand{\mStress}{\sigma}
\newcommand{\Mmaterialtangent}{C}
\newcommand{\Biot}{\mathrm{B}}

\newcommand{\mdisplacement}{u}
\newcommand{\mposition}{x}

\newcommand{\ROMbasis}{\phi}
\newcommand{\snapshot}{x}
\newcommand{\mwork}{w}

\newcommand{\mstiffness}{\Matr{k}}
\newcommand{\mintforce}{\matr[\mathrm{int}]{\nodal{f}}}

\newcommand{\mResidual}{\matr{r}}
\journal{Computer Methods in Applied Mechanics and Engineering }

\begin{document}
\begin{frontmatter}
\title{A Monolithic Hyper ROM FE\textsuperscript{2} Method with Clustered Training at Finite Deformations}
\author[Freiberg]{Nils Lange\footnote[2]{Corresponding author. \textit{E-mail address:} \href{Nils.Lange@imfd.tu-freiberg.de}{Nils.Lange@imfd.tu-freiberg.de} }}
\author[Freiberg,Cottbus]{Geralf Hütter}
\author[Freiberg]{Bjoern Kiefer}
\affiliation[Freiberg]{organization={TU Bergakademie Freiberg, Institute of Mechanics and Fluid Dynamics},%Department and Organization
            addressline={Lampadiusstr.~4}, 
            city={Freiberg},
            postcode={09596}, 
            country={Germany}}
\affiliation[Cottbus]{organization={BTU Cottbus-Senftenberg, Institute of Civil Engineering},%Department and Organization
            addressline={Konrad-Wachsmann-Allee 2}, 
            city={Cottbus},
            postcode={03046}, 
            country={Germany}}
\begin{abstract}
The usage of numerical homogenization to obtain structure-property relations using the finite element method at both the micro and macroscale has gained much interest in the research community. However the computational cost of this so called FE\textsuperscript{2} method is so high that algorithmic modifications and reduction methods are essential. Currently the authors proposed a monolithic algorithm. Now this algorithm is combined with ROM and ECM hyper integration, applied at finite deformations and complemented by a clustered training strategy, which lowers the training effort and the number of necessary ROM modes immensely. The applied methods are modularly combinable as aimed in finite element approaches. An implementation in terms of an extension for the already established MonolithFE\textsuperscript{2} code is provided. Numerical examples show the efficiency and accuracy of the monolithic hyper ROM FE\textsuperscript{2} method and the advantages of the clustered training strategy. Online times of below $1\%$ of the conventional FE\textsuperscript{2} method could be gained. In addition the training stage requires around $3\%$  of that time, meaning that no extremely expensive offline stage is necessary as in many Neural Network approaches, which only pay off when a lot of online simulations will be conducted.
\end{abstract}

\begin{graphicalabstract}
\begin{figure}
\centering
\includegraphics[width=\textwidth]{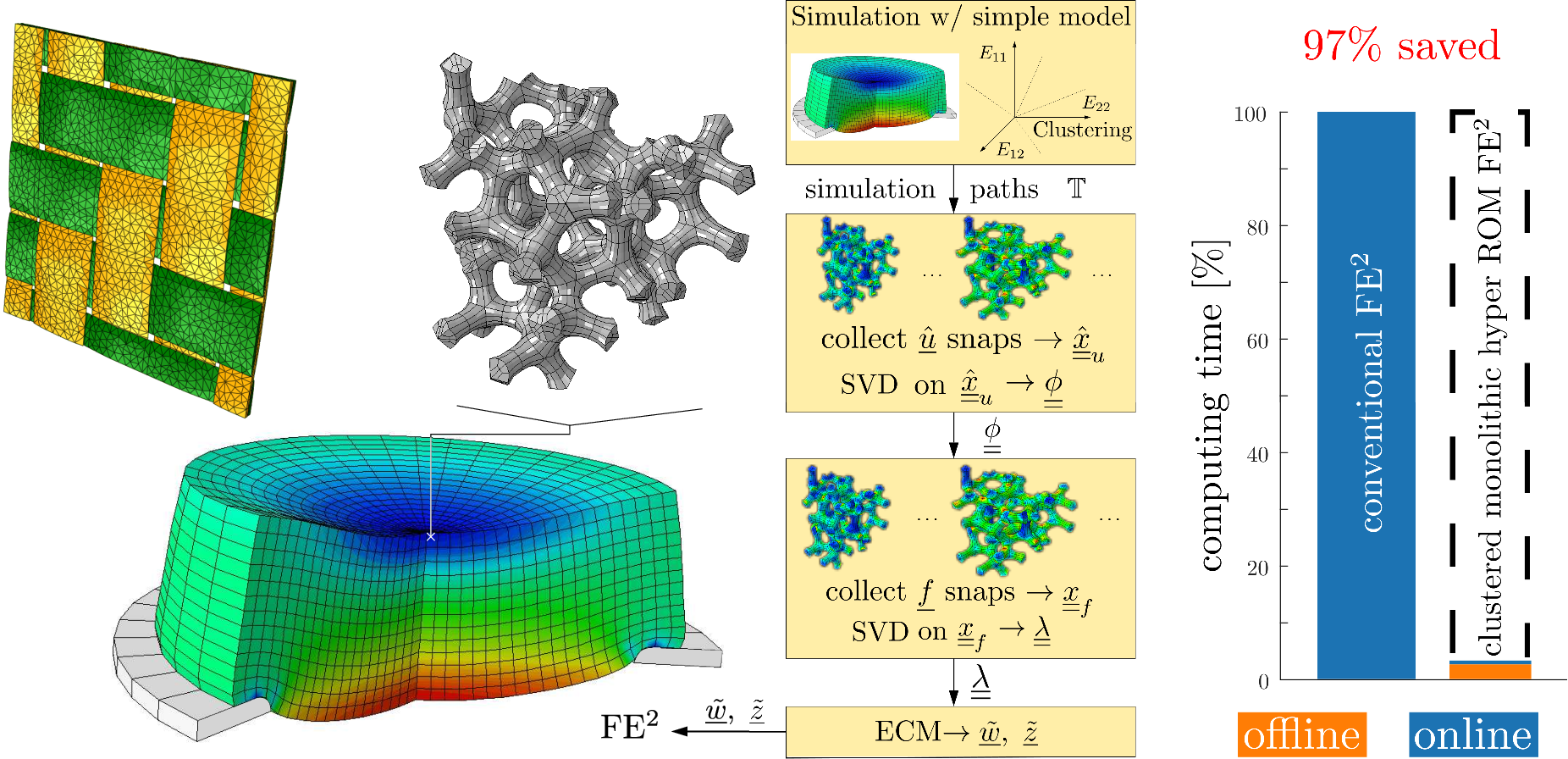}
\end{figure}
\end{graphicalabstract}

\begin{keyword}
Homogenization\sep Multiscale \sep FE\textsuperscript{2}\sep ROM\sep Hyper Integration \sep Monolithic Scheme

\end{keyword}

\end{frontmatter}

%\begin{small}
%\tableofcontents   %\end{small}

\section{Introduction}
Nearly all materials relevant for engineering show characteristic inhomogeneities at a certain length scale. Aiming to describe the materials behavior in an adequate accuracy, it usually is sufficient to use homogenized material laws working with a strain definition as observable state variable and a set of internal state variables with evolution laws for the associated forces, in order to describe all macroscopically observed effects, while ensuring thermodynamic consistency. For certain very complex microstructures e.g.~composites or foams, homogenized materials may be hard to find, but usually the behavior of the individual components is known. Resolving the microstructure fully in a structural model is usually unfavorable, since firstly the microstructure is often not known everywhere in the sense of actual data e.g.~from CT examinations and secondly the computational effort becomes tremendously high for large scale separations.

The necessity of solving these problems often leads to the common way of multiscale modeling using the assumption of a clear scale separation, which may be seen as a model reduction step. It is hypothesized that a representative volume element (RVE) can be found, including all essential microstructural defects and inclusions and showing statistically the same mechanical behavior as the actual microstructure, which describes the behavior of the point of a structure using the principle of local action.

The idea of examining the scale transition between a heterogeneous material and a structural level are of fundamental interest. This can be seen from the fact that basic concepts already go back to the second half of the 19\textsuperscript{th} century. From this time originates the Voigt bound, later the Reuss, Hill, Kr\"oner etc.~bounds were proposed, an overview can be found in \citep{Hirsekorn_1990}. An important fundamental consideration concerning the scale transition was the condition of macro homogeneity \citep{Hill1972}. Later more distinct analytical homogenization approaches like the Mori-Tanaka method \citep{Mori_Tanaka_1973} were proposed. The computational resources in this time period were too restricted to even consider concurrent multiscale concepts. Due to the exponential growth of computer power according to Moores law, F.~Feyel was able to present an implementation in the late 1990s \citep{feyel_multiscaleFE2_1999}, which linked a coarse and fine scale, whereat on both levels the finite element method (FEM) was used, for which he coined the term \enquote{FE\textsuperscript{2} method}. This method developed out of an engineering research environment, whereas more mathematically oriented researches came to a similar approach called \enquote{FE-HMM}, cf.~a comparison of the approaches and a nonlinear version of the FE-HMM in \citep{Eidel_2018,Eidel_2019}. Other multiscale methods that leave behind the continuum assumption e.g.~peridynamics that are also frequently used, are outside of the scope of this article. In view of the enormous amount of research activity in the field of multiscale modeling, overview articles like \citep{kanoute_multiscale_2009,Geers2010,Schroeder2014,Geers_multiscale_2017,Matous2017,Raju2021,Fish_2021} are necessary to classify and compare the different concepts.

Recently more and more research groups with different backgrounds pour into the scene, since it has been observed that the numerical multiscale concept using homogenization is not only suited for mechanical problems, but also for electro-magnetic, chemical, biological issues and their interrelationship. Regarding the placement of the approaches into Gartner's hype cycle, \citet{Fish_2021} rate the current state as a phase of \enquote{disillusionment}. The expectations starting from the early innovation triggers rose to a hype of thinking \enquote{with these new simulation methods one can simulate \emph{everything}}. Now a state of depression has arrived in the collective thinking, because it is realized that the computational resources are often too restricted for even the smallest realistic problems. This stems from the fact that if using the finite element method on both levels, then behind every macroscopic integration point a very expensive microscopic FE model has to be computed. An incomplete list of the common endeavors to reduce the computational effort is:
\begin{enumerate}
     \setlength{\itemsep}{1pt}
    \item Macroscopic elements with reduced integration
    \item Approximation of the macro tangent e.g.~Quasi-Newton algorithm
    \item Massive parallelization in the macroscopic element loop
    \item Usage of a monolithic (\enquote{true}) Newton scheme instead of a staggered one \label{enum:monolitic}
    \item Fast-Fourier Transformation (FFT) solver on the micro level \label{enum:FFT}
    \item Description of the macroscopic material behavior via Neural Networks trained with micro FE simulations\label{enum:NeuralNetwork}
    \item Application of Reduced Order Modeling (ROM) together with hyper integration at the micro level \label{enum:hyperROM}
\end{enumerate}

While the first \ref{enum:FFT} ideas represent only algorithmic variations of the original method, approaches \ref{enum:NeuralNetwork} and \ref{enum:hyperROM} have an approximative character. Neural Networks (NN) are quite common in many recent applications. When the material response is nonlinear elastic, they are capable of describing a material very accurately and efficiently. When the material response is irreversible, much more complicated considerations must be taken into account. There are e.g.~hybrid models \citep{Settgast2020,Jang2020,Qu2021,Malik2022,Kalina2022}, or so-called deep material networks \citep{Liu2019,Gajek2020} to ensure thermodynamic consistency. A review article classifying these methods can be found in \citep{Liu2021}, in which it is stressed that certain constraints must be present to ensure that physically plausible and interpretable results will be obtained. These considerations however must be taken into account for every class of material and are hardly to be generalized. The idea \ref{enum:hyperROM}, which was to the best knowledge of the authors first applied by \citet{Hernandez2014} in the nonlinear/irreversible multiscale context and picked up by various other authors e.g.~\citet{Caicedo2019,Rocha2019,Soldner2016}, however directly ensures a physical material description especially under irreversible behavior. Thereby the properties of the original FE\textsuperscript{2} model are being kept without the need for special considerations, which is why this article only focuses on this approach. It is worth mentioning that a similar method was also proposed for the FFT formulation, whereby in complete analogy a reduced set of Fourier Modes is constructed and a lower set of integration points is chosen \citep{Waimann2022}. Other sophisticated methods relying also on the reduced basis approach but moving further away from the FE\textsuperscript{2} idea towards hybrid approaches e.g.~the works of \citet{Fritzen2013,Fritzen2019}, \citet{Sharba2023} and the usage of multiscale ROM methods outside the sole mechanical context e.g.~in the analysis of heterogeneous poroelastic media \citep{Jaenicke2016}, are outside of the scope of this research.

This article has six goals. At first an efficient way of handling finite deformations in the context of the ROM FE\textsuperscript{2} method is presented. Secondly the incorporation of the monolithic algorithm into the ROM method is shown. Thirdly the original proposal of the Empirical Cubature Method (ECM) by Hern\'{a}ndez et al.~is revisited and applied. Fourthly a clustered training strategy is introduced in order to reduce the training cost and fifthly an implementation for the commercial software Abaqus is shown. Finally numerical examples are provided, proving the computational efficiency and accuracy of the presented approaches.

\section{Multiscale Modeling}
The main idea of multiscale modeling lies in the fact that often specific distinct scales with a large separation of their characteristic lengths exist, that influence each other, but whose local phenomena are to be modeled at the correct scale to get a meaningful effect description. Scales are e.g.~the interaction between different galaxies, sun systems in a galaxy, planets in a sun system, movement and deformation of body's on a planet, the deformation of the meso structure of the material of a body, the movement of dislocations and atoms of a material, the interaction of elementary particles and so on. Usually it is not constructive to incorporate the local behavior of a specific scale in a higher scale, e.g.~the attempt of describing the movement of galaxies by including the behavior of all the individual atoms in the modeling approach would be unreasonable, even though it might be possible in general. The usual approach in physics is to make an experiment (sometimes replaced by simulations) on the desired scale to get an effective behavior, even if it is known that this effective behavior depends on the structure and properties of the lower scales.

\begin{figure}[!h]
	\centering
	\includegraphics[width=0.58\textwidth]{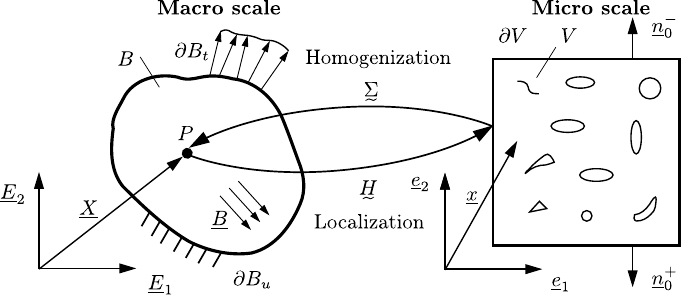}
	\caption{Mechanical two-level multiscale model: Localization of the locally uniform macroscopic deformation state, Homogenization of the inhomogeneous microscopic stress response}
	\label{fig:multiscale_scheme}
\end{figure}

In material modeling the effective properties can be very complicated, especially when the material behavior is strongly irreversible and experiments might be very expensive and complicated. Since engineering materials are artificial, its microstructure and the behavior of the individual components are sometimes well known and understood. Then the at the structural scale well approved Cauchy continuum might be just applied at the microscale (depending on the size also called meso scale). The multiscale problem is depicted in Figure \ref{fig:multiscale_scheme}. It is assumed that the macroscopic values of mechanical power, 1st Piola Kirchhoff stress $\undertilde{\Sigma}^{\mathrm{PK}}$ and displacement gradient $\undertilde{H}$ can just be defined as the volume average of their microscopic counterparts, whereat the first is known as \enquote{Hill-Mandel} or \enquote{macro homogeneity} condition:
\begin{equation}
   J\,\undertilde{\Sigma}:\undertilde{D}=\undertilde{\Sigma}^{\mathrm{PK}}:\undertilde{\dot{H}}=\langle\undertilde{\sigma}^{\mathrm{PK}}:\undertilde{\dot{h}}{\rangle}_0
\end{equation}
\begin{equation}\label{eqn:homogenization_Stress}
    \undertilde{\Sigma}^{\mathrm{PK}}=\langle\undertilde{\sigma}^{\mathrm{PK}}{\rangle}_0=\dfrac{1}{V_0}\int_{\partial V_0}\underline{t}\otimes\underline{x}_0\ \mathrm{d}\,S
\end{equation}
\begin{equation}\label{eqn:homogenization_Strain}
    \undertilde{H}=\underline{U}\otimes\nabla_{\underline{X}_0}=\langle\undertilde{h}{\rangle}_0=\dfrac{1}{V_0}\int_{\partial V_0}\underline{u}\otimes\underline{n}_0\ \mathrm{d}\,S
\end{equation}
A RVE must be identified, being representative for the surrounding of a point $P$. The scale bridging is established in terms of Homogenization by the definition of the stress in form of the surface integral as shown in equation \eqref{eqn:homogenization_Stress}. The boundary value problem at the macroscale reads:
\begin{equation}
    \mathrm{div}(\undertilde{\Sigma})+\underline{B}=\underline{0}\ \mathrm{in}\ B\ ,\ \ \ \underline{U}=\underline{\bar{U}}\ \ \mathrm{on}\ \ \partial B_u\ ,\ \ \ \underline{T}=\underline{\bar{T}}\ \ \mathrm{on}\ \ \partial B_t\  
\end{equation}
The Localization step is introduced, by assuring that equation \eqref{eqn:homogenization_Strain} is fulfilled a priori through suitable boundary conditions, whereby usually periodic boundary conditions are chosen, which assume a periodic RVE continuation by identifying homologous points $()^{{\small{+}}}$ and $()^{{\small{-}}}$. All together the microscopic boundary value problem reads:
\begin{equation}\label{eqn:FEM_micro}
    \mathrm{div}(\undertilde{\sigma})=\underline{0}\ \mathrm{in}\ V\ ,\ \ \ \ \underline{u}^{{\small{+}}}-\underline{u}^{{\small{-}}}=\undertilde{H}\cdot\left[\underline{x}^{{\small{+}}}_0-\underline{x}_0^{{\small{-}}}\right]\ \ \mathrm{on}\ \ \partial V\ \ \mathrm{with}\ \ \underline{n}^{{\small{-}}}_0=-\underline{n}^{{\small{+}}}_0
\end{equation}
The partial differential equations of the involved problems must be solved numerically, usually by means of the FEM. The discretized equilibrium equations of the macro problem in residual form read:
\begin{equation}
    \underline{\hat{R}}\,(\underline{\hat{U}},\underline{\hat{u}})=\underbrace{\sum_{\alpha=1}^{i}\underline{\underline{B}}^\mathrm{T}_{\,\alpha}\cdot\underline{\Sigma}_{\,\alpha}\,(\underline{H}_{\,\alpha}\!=\!\underline{\underline{B}}_{\,\alpha}\cdot\underline{\hat{U}},\underline{\hat{u}}_{\,\alpha})}_{\textstyle\underline{\hat{F}}_\mathrm{\,int}}-\underline{\hat{F}}_\mathrm{\,ext}=\underline{0}
\end{equation}
Thereby all measures with $(\,\hat{}\,)$ denote nodal values. The numerical integration interrogates the stress response of an attached RVE at all integration points $\alpha$ from $1$ to $i$. The projection onto the global B-Matrix yields the vector of internal forces. At the microscale the discretized FEM equations read for one macroscopic integration point $\alpha$:
\begin{equation}
    \underline{\hat{u}}^{{\small{-}}}_{\,\alpha}\,(\underline{\hat{u}}_{\,\alpha},\underline{H}_{\,\alpha})=\underline{\underline{A}}\cdot\begin{bmatrix}\underline{\hat{u}}_{\,\alpha}\\ \underline{H}_{\,\alpha}\end{bmatrix},\ \underline{\hat{u}}_{\,\alpha}=\begin{bmatrix}\underline{\hat{u}}^{{\small{+}}}_{\,\alpha}\\ \underline{\hat{u}}^{\mathrm{inner}}_{\,\alpha}\end{bmatrix},\ \underline{\hat{r}}\,(\underline{H}_{\,\alpha}(\underline{\hat{U}}),\underline{\hat{u}}_{\,\alpha})=\underline{\hat{f}}_{\,\mathrm{int}\,\alpha}+\underline{\underline{A}}^\mathrm{T}\cdot\underline{\hat{f}}_{\,\mathrm{int}\,\alpha}^{{\,\small{-}}}=\begin{bmatrix}\underline{0}\\\underline{\Sigma}^\mathrm{PK}_{\,\alpha}\cdot V_0\end{bmatrix}
\end{equation}
Therein, a matrix $\underline{\underline{A}}$ is introduced giving the linear relation \eqref{eqn:FEM_micro} between the displacement of the \enquote{$^-$} nodes and the displacement of the \enquote{$^+$} nodes and the macro displacement gradient. The superscript $^{\mathrm{inner}}$ corresponds to the nodes which are not at the boundary. The Galerkin projection for the \enquote{$^-$} internal force with $\underline{\underline{A}}$ then reflects the fact that $\underline{t}^{{\small{-}}}=-\underline{t}^{{\small{+}}}$ and that the force response to the macro displacement gradient will be the macro stress.

\section{Solution Strategy}
\subsection{General aspects} 
\begin{figure}[!h]
  \centering
  \subfloat[][Computational efficiency vs.~flexibility for different methods.\label{fig:comparison_computational_methods}]{\includegraphics[width=0.6\linewidth]{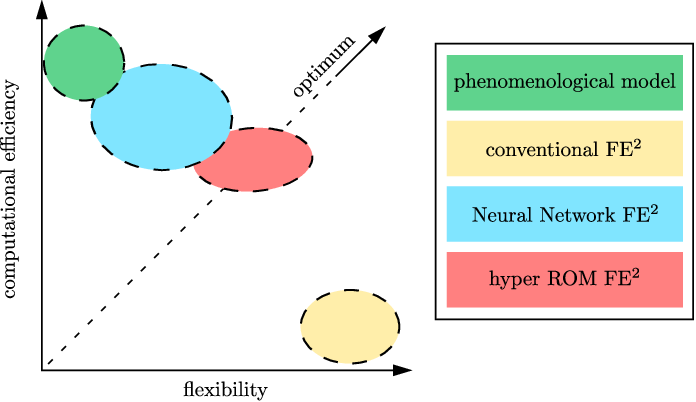}}
  \hspace{0.04\linewidth}
  \subfloat[][Composition of modeling effort.\label{fig:composition_computational_effort}]{\includegraphics[width=0.35\linewidth]{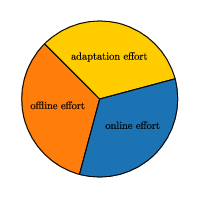}}
  \caption{General classification of formulations for materials with characteristic microstructure.}
  \label{fig:Computational_methods}
\end{figure}

In the past decades, many approaches  for multiscale material modeling have emerged. For now, we want to compare them in a general sense, namely in terms of flexibility and computational efficiency, assuming all methods to be as accurate as needed for engineering purposes. The optimal material formulation would be flexible to account for any type of material, while needing as few computational resources as possible, cf.~Figure~\ref{fig:Computational_methods}(\subref{fig:comparison_computational_methods}). As illustrated in Figure~\ref{fig:Computational_methods}(\subref{fig:composition_computational_effort}), the effort in modeling can be split into online costs, which include all needed resources for the actual simulations, offline costs, which include possible training data generation and evaluation as well as adaptation effort, which comprises all necessary steps to develop or modify a method to account for a new specific material behavior.

Classical phenomenological material descriptions, that can also be micromechanically motivated, such as the Gurson-Tvergaard-Needleman (GTN) model for porous metals \cite{Tvergaard1984} or the Deshpande-Fleck model for foams \cite{Deshpande2000}, are often very efficient, but limited in the sense that they only capture the characteristic effects of specific material classes. In contrast, the conventional FE\textsuperscript{2} method, see  Feyel \cite{feyel_multiscaleFE2_1999}, in this context often referred to as High Fidelity (HF) model, is theoretically able to account for any type of material, but is so computationally inefficient that it may never leave the academic context. In order to speed up FE\textsuperscript{2} computations, Neuronal Network based material formulations have been developed, which rely on RVE-FE training data as input. They have been shown to save a lot of computational costs in the online phase \citep{Rocha2020}. On the downside, they require a lot of training data, as Neuronal Networks do not extrapolate well. Moreover, since they only indirectly represent solutions to the original PDEs of the underlying problem, careful considerations must be taken into account to obtain physically consistent responses, and therefore also additional adaptation effort \citep{Liu2021}. This, in  turn, may then result in formulations that may still be able to account for a broader spectrum of problems than classical material descriptions, but are much more restricted than the conventional FE\textsuperscript{2} method. The hyper ROM FE\textsuperscript{2}, on the other hand, still solves the associated system of PDEs and extra- and interpolates well in untrained regions, so that it generally needs less training data and no adaptation effort is necessary to account for specific material behavior. The online costs, however, are usually higher than those of Neuronal Network based models \citep{Rocha2020}.

\subsection{Reduced Order Modeling}\label{sec:ROM}
\subsubsection{Basic concepts}
Reduced Order Modeling is a class of techniques used quite frequently in a lot of different fields of application. Its main idea lies in the observation that many problems containing a lot of information can usually be described with very good accuracy by so called lower order approximations, because the degree of correlation is high. In particular, considering the recent trend of gathering more and more data, the problem often arises of having to evaluate this data in reasonable time or to store the information in a compressed way. The roots go back to the 19\textsuperscript{th} century, whereby depending on the field of application different names for the same/similar technique are established for decomposing large amounts of data with a high amount of coherence, including Proper Orthogonal Decomposition (POD), Singular Value Decomposition (SVD), Principal Component Analysis (PCA), Karhunen-Loeve Transforma (KLT), Hotelling Transform among others. A broad overview can be found in \citep{Bebendorf2014}.

Regarding a (discrete) mechanical problem, where the value of the nodal displacements $\matr{\nodal{\mdisplacement}}(t)\in\mathbb{R}^{\nodal{n}}$ of an RVE at time $t$ are being sought, this kind of approach leads to linear combinations of orthonormal basis vectors $\matr{\ROMbasis}_i\in \mathbb{R}^{\nodal{n}}, \ i=1,\dots,\amplitudal{n}$ with $\amplitudal{n}\ll \nodal{n}$, yielding
\begin{equation}\label{ROM_assumption}
    \matr{\nodal{\mdisplacement}}(t)\approx\Matr{\ROMbasis}\cdot\matr{\amplitudal{\mdisplacement}}(t)\ ,
\end{equation}
\noindent wherein $\matr{\amplitudal{\mdisplacement}}\in\mathbb{R}^{\amplitudal{n}}$ are the amplitudes of the corresponding basis vectors. Here, and in the following, a superimposed tilde ($\tilde{\phantom{a}}$) marks the reduced representation of a quantity. The basis vectors $\matr{\ROMbasis}_i$ can be interpreted as deformation modes, of which only the most dominant ones are being considered. The modes can be found through a SVD applied to a collection of displacement snapshots from RVE simulations using representative loading paths, yielding the matrix decomposition as shown in equation \eqref{SVD_displacement}. In doing so, an elastic-inelastic split, as proposed by \citet{Hernandez2014} should be used, which assures that the elastic range is recovered exactly:
\begin{equation}\label{SVD_displacement}
    \Matr[\mathrm{el}]{\nodal{\snapshot}}_\mdisplacement=\Matr[\mathrm{el}]{\ROMbasis}\cdot\Matr[\mathrm{el}]{\sigma}_\mdisplacement\cdot\underline{\underline{v}}^{\mathrm{el\,T}}_u\ ,\ \ \ \ \Matr[\mathrm{inel\,mod}]{\nodal{\snapshot}}_\mdisplacement=\Matr[\mathrm{inel}]{\ROMbasis}\cdot\Matr[\mathrm{inel}]{\sigma}_\mdisplacement\cdot\underline{\underline{v}}^{\mathrm{inel\,T}}_u\ ,\ \ \ \ \Matr{\ROMbasis}=\left[\Matr[\mathrm{el}]{\ROMbasis}\ \Matr[\mathrm{inel}]{\ROMbasis}\right]\ .
\end{equation}
\noindent Here, $\Matr{\nodal{\snapshot}}_\mdisplacement$ contains the value of the nodal displacements at selected points of time from each training load case, i.e.
\begin{equation}\label{snapshots_displacement}
    \Matr[\mathrm{el}]{\nodal{\snapshot}}_\mdisplacement=\left[\matr[\mathrm{el}\,1]{\nodal{\mdisplacement}},\,\matr[\mathrm{el}\,2]{\nodal{\mdisplacement}},...\,\right]\ ,\ \ \ \ \Matr[\mathrm{inel}]{\nodal{\snapshot}}_\mdisplacement=\left[\matr[\mathrm{inel}\,1]{\nodal{\mdisplacement}},\,\matr[\mathrm{inel}\,2]{\nodal{\mdisplacement}},...\,\right]\ ,\ \ \ \ \Matr[\mathrm{inel\,mod}]{\nodal{\snapshot}}_\mdisplacement=\Matr[\mathrm{inel}]{\nodal{\snapshot}}_\mdisplacement-\Matr[\mathrm{el}]{\ROMbasis}\cdot\transp{\Matr[\mathrm{el}]{\ROMbasis}}\cdot\Matr[\mathrm{inel}]{\nodal{\snapshot}}_\mdisplacement\ .
\end{equation}
\noindent From equation \eqref{ROM_assumption}, it follows for the nodal coefficients of the test field in the Galerkin Method $\delta\vect{\mdisplacement}(\vect{\mposition})$
\begin{equation}
    \delta\matr{\nodal{\mdisplacement}}(t)\approx\Matr{\ROMbasis}\cdot\delta\matr{\amplitudal{\mdisplacement}}(t)\ ,
\end{equation}
\noindent since the same ansatz is being made for the actual and test displacement field. Considering the virtual work of an RVE, assuming a quasi-static problem with a divergence-free stress field, this leads to
\begin{equation}\label{ROMintforce}
    \delta\mwork=0=\delta\transp{\matr{\nodal{\mdisplacement}}}\cdot\int_V\transp{\Matr{B}}\cdot\matr{\mStress}\,\mathrm{d}V=\delta\transp{\matr{\nodal{\mdisplacement}}}\cdot\mintforce=\delta\transp{\matr{\nodal{\mdisplacement}}}\cdot\nodal{\mResidual}=\delta\transp{\matr{\amplitudal{\mdisplacement}}}\cdot\transp{\Matr{\ROMbasis}}\cdot\nodal{\mResidual}=\delta\transp{\matr{\amplitudal{\mdisplacement}}}\cdot\amplitudal{\mResidual}\ .
\end{equation}
\noindent The linearization of equation \eqref{ROMintforce} yields
\begin{equation}\label{ROMstiffness}
    \dfracderiv{\amplitudal{\mResidual}}{\matr{\amplitudal{\mdisplacement}}}=\dfracderiv{\amplitudal{\mResidual}}{\nodal{\mResidual}}\cdot\dfracderiv{\nodal{\mResidual}}{\matr{\nodal{\mdisplacement}}}\cdot\dfracderiv{\matr{\nodal{\mdisplacement}}}{\matr{\amplitudal{\mdisplacement}}}=\transp{\Matr{\ROMbasis}}\cdot\nodal{\mstiffness}\cdot\Matr{\ROMbasis}=\amplitudal{\mstiffness}\ ,
\end{equation}
whereby $\amplitudal{\mstiffness}$ becomes a dense matrix.

\subsubsection{Finite deformations}
When considering a large strain theory, a rigid body rotation is superimposed onto an RVE in addition to the deformation. Then either the rotations must be trained, by selecting the deformation gradient as deformation measure, as proposed by \citet{Caicedo2019}, or training directions, that do not contain macroscopic rigid body motion are chosen. The latter  requires the basis system to be rotated along in the online stage, since the ROM basis is then defined in an unrotated configuration. 

This, however  may be impractical and time consuming in an actual implementation. Instead, it seems promising to only impose the stretch onto the RVE and obtain the (symmetric) Biot stress in response. Working with these symmetric stress/strain measures additionally has the advantage of lowering the computational effort in computing the (macro-)stress and condensing the tangent (in 3D only 6 instead of 9 independent variables). This was also observed by \citet{Kochmann_2018}. No macroscopic rotations have to be trained in this approach. The additional benefit is that this allows a code to be quickly adapted, with little effort to make it compatible to both small and large deformation theory. 

To this end, the polar decomposition of the deformation gradient is considered first
\begin{equation}\label{eqn:ploar_decomp}
    \undertilde{F}=\undertilde{H}+\undertilde{I}=\undertilde{R}\cdot\undertilde{U}\ ,\ \undertilde{L}=\undertilde{\dot{F}}\cdot\undertilde{F}^{-1}=\undertilde{D}+\undertilde{W}=\underbrace{\undertilde{R}\cdot\mathrm{sym}(\undertilde{\dot{U}}\cdot\undertilde{U}^{-1})\cdot\undertilde{R}^\mathrm{T}}_{\textstyle\undertilde{D}}+\underbrace{\undertilde{\dot{R}}\cdot\undertilde{R}^\mathrm{T}+\undertilde{R}\cdot\mathrm{skw}(\undertilde{\dot{U}}\cdot\undertilde{U}^{-1})\cdot\undertilde{R}^\mathrm{T}}_{\textstyle\undertilde{W}},
\end{equation}
\noindent where $\tens{\MRotation}$ is the macro rotation, $\tens{\MStretch}$ the right stretch, $\tens{D}$ the rate of deformation and $\tens{W}$ the spin tensor. Considering the mechanical power, it can be shown that $\undertilde{\dot{U}}$ is work conjugate to the (symmetric) Biot stress $\tens[\Biot]{\MStress}$:
\begin{equation}\label{workconjugacy}
	%\tens[\Piola]{\MStress}:\tens{\MDefgrad}=\tens[\Piola]{\MStress}:[\tens{\MRotation}\cdot\tens{\MStretch}]=[\transp{\tens{\MRotation}}\cdot\tens[\Piola]{\MStress}]:\tens{\MStretch}\overset{\tens{\MStretch}=\transp{\tens{\MStretch}}}{=}\dfrac{1}{2}[\transp{\tens{\MRotation}}\cdot\tens[\Piola]{\MStress}+\transp{\tens[\Piola]{\MStress}}\cdot\tens{\MRotation}]:\tens{\MStretch}=\tens[\Biot]{\MStress}:\tens{\MStretch}
	J\,\undertilde{\Sigma}:\undertilde{D}=J\,\undertilde{\Sigma}:\undertilde{R}\cdot\mathrm{sym}(\undertilde{\dot{U}}\cdot\undertilde{U}^{-1})\cdot\undertilde{R}^\mathrm{T}=\mathrm{sym}(\undertilde{R}^\mathrm{T}\cdot J\,\undertilde{\Sigma}\cdot\undertilde{F}^{-\mathrm{T}}):\undertilde{\dot{U}}=\tens[\Biot]{\MStress}:\undertilde{\dot{U}}\ .
\end{equation}
Since a lot of commercially available FE codes, such Abaqus and Ansys, work with the true (Cauchy) stress, the output in terms of stress and tangent must be transformed by using push-forward operations. Considering equation \eqref{workconjugacy}, the Cauchy stress can then be written as a function of the Biot stress, namely
\begin{equation}\label{biotstress}
	\tens[\Biot]{\MStress}=\dfrac{\MJacobian}{2}[\transp{\tens{\MRotation}}\cdot\tens{\MStress}\cdot\transpinv{\tens{\MDefgrad}}+\invers{\tens{\MDefgrad}}\cdot\tens{\MStress}\cdot\tens{\MRotation}]\ \leftrightarrow\  \undertilde{\Sigma}=\dfrac{2}{\MJacobian}\invers{[\indexnot{\MRotation}{ik}\invers{\indexnot{\MDefgrad}{lj}}+\invers{\indexnot{\MDefgrad}{ki}}\indexnot{\MRotation}{jl}]}\indexnot[\Biot]{\MStress}{kl}\, \underline{E}_i\otimes\underline{E}_j
\end{equation}
It is further shown in \ref{sec:appendix_biot_tangent} how to compute the necessary algorithmically consistent tangent, as required in the Updated Lagrangian approach context.

\subsubsection{Monolithic solution scheme}
\begin{figure}[!h]
	\centering
	\includegraphics[width=0.75\textwidth]{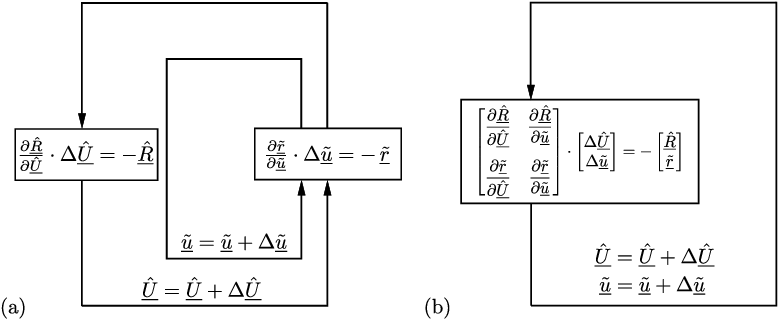}
	\caption{(a) Staggered algorithm, (b) monolithic algorithm, see also \citep{Lange2021_monolithic}}
	\label{fig:monolithic_staggered}
\end{figure}

The numerical treatment in solving equation \eqref{ROMintforce} follows in analogy to the full problem. Regarding the monolithic scheme previously by the authors proposed \citep{Lange2021_monolithic}, see Figure \ref{fig:monolithic_staggered}(b), the only change is in using amplitudal quantities instead of nodal ones at the microscale.

For the increment of the microscopic nodal displacements this yields
\begin{equation}
    \Delta\vect{\amplitudal{\mdisplacement}}=\,-\invers{\left[\,\transp{\Matr{\ROMbasis}}\cdot\nodal{\mstiffness}\cdot\Matr{\ROMbasis}\,\right]}\cdot\left[\,\transp{\Matr{\ROMbasis}}\cdot\nodal{\matr{\mResidual}}+\transp{\Matr{\ROMbasis}}\cdot\dfracderivpartial{\nodal{\mResidual}}{\matr{\MStretch}}\cdot\Delta\matr{\MStretch}\,\right]=-\,\invers{\amplitudal{\mstiffness}}\cdot\left[\,\amplitudal{\matr{\mResidual}}+\dfracderivpartial{\amplitudal{\mResidual}}{\matr{\MStretch}}\cdot\Delta\matr{\MStretch}\,\right]\ ,
\end{equation}
wherein $\Delta\matr{\MStretch}$ is the increment of the macroscopic stretch over one Newton step. The algorithmically consistent macroscopic (Biot) stress can be computed as
\begin{equation}
    \matr[\Biot]{\MStress}_{\mathrm{alg}}=\matr[\Biot]{\MStress}-\dfracderivpartial{\matr[\Biot]{\MStress}}{\matr{\nodal{\mdisplacement}}}\cdot\Matr{\ROMbasis}\cdot\invers{\left[\,\transp{\Matr{\ROMbasis}}\cdot\nodal{\mstiffness}\cdot\Matr{\ROMbasis}\,\right]}\cdot\transp{\Matr{\ROMbasis}}\cdot\nodal{\matr{\mResidual}}=\matr[\Biot]{\MStress}-\dfracderivpartial{\matr[\Biot]{\MStress}}{\matr{\amplitudal{\mdisplacement}}}\cdot\invers{\amplitudal{\mstiffness}}\cdot\amplitudal{\mResidual}\ .
\end{equation} The tangent stiffness reads in terms of amplitudal quantities
\begin{equation}\label{Biot_Tangent}
    \Matr[\Biot]{\Mmaterialtangent}=\dfracderivpartial{\matr[\Biot]{\MStress}}{\dot{\matr{\MStretch}}}-\dfracderivpartial{\matr[\Biot]{\MStress}}{\matr{\nodal{\mdisplacement}}}\cdot\Matr{\ROMbasis}\cdot\invers{\left[\,\transp{\Matr{\ROMbasis}}\cdot\nodal{\mstiffness}\cdot\Matr{\ROMbasis}\,\right]}\cdot\transp{\Matr{\ROMbasis}}\cdot\dfracderivpartial{\nodal{{\mResidual}}}{\dot{\matr{\MStretch}}}=\dfracderivpartial{\matr[\Biot]{\MStress}}{\dot{\matr{\MStretch}}}-\dfracderivpartial{\matr[\Biot]{\MStress}}{\matr{\amplitudal{\mdisplacement}}}\cdot\invers{\amplitudal{\mstiffness}}\cdot\dfracderivpartial{\amplitudal{{\mResidual}}}{\dot{\matr{\MStretch}}}\ .
\end{equation}

\subsection{Hyper integration}
The main contribution of the approach discussed in Section~\ref{sec:ROM} towards computational efficiency was that only a small number of $\amplitudal{n}$ equations have to be solved instead  of the large number of $\nodal{n}$ DOF of the full system. Nevertheless, evaluating the integral expression in equation \eqref{ROMintforce}, as well as computing the macroscopic stresses is still very expensive.  Moreover, the projection with $\Matr{\ROMbasis}$ actually introduces additional computational effort.
%The result of the ROM Galerkin projection is, that the internal force of one node contributes to the internal force of all modes. Therefore also all integration points contribute to the force of all modes.
%In a linear elastic problem, only 6 (2D:3) modes are needed to describe all possible deformation states. Per mode in general now a linear combination of the internal force of at most 6 (2D:3) integration points, whose force vectors are linearly independent, would be sufficient to integrate the problem exactly, making it in total maximal 36 (2D:9) integration points to integrate the force of all possible states exactly. It is thereby not important which integration points are chosen. It is enough to find \textit{any} set of integration points whose internal force contribution can integrate the modal forces by their linear combination.

There are mainly two methods described in the literature for integrating the reduced order problem, namely the Empirical Interpolation Methods (EIM) and the Empirical Cubature Methods (ECM), cf.~\citep{Hernandez2014,Tuijl2018}. A comparison of their numerical performances in the multiscale context can, for instance, be found in \citep{Tuijl2018,Soldner2017}. Therein, \citet{Tuijl2018} conclude that even though EIM methods lead to higher speedups, they tend to produce more interpolation errors. It is also less straightforward to construct accurate snapshots spaces in comparison to the ECM methods, which are therefore more promising in that regard. In the scope of this work, we will therefore only focus on the ECM approach.

In the context of multiscale modeling, \citet{Hernandez2014} were the first to comprehensively adapt the ECM to the needs of the FE\textsuperscript{2} method, but the idea traces back to the work of \citet{An2009}. An overview of further roots can be found in \citep{Hernandez2016}. Subsequently, the algorithm was applied by different authors, e.g.~\cite{Rocha2018,Tuijl2018}, and also further developed. This resulted, for instance, in the  Reduced Energy Optimal Qubature (REOQ) method \cite{Caicedo2019,Raschi2021}, which relies on the integration of the free energy. \citet{Hernandez2020} later applied his method  to general periodic structures, which yielded an effective tool for weak scale separations, but cannot be reasonably applied for large scale separations. An interesting work of \citet{Rocha2019} presented an offline-free adaptive ROM/ECM implementation, which does not require a training stage, but yields less online speedup than the original method of Hern\'{a}ndez et al.

\noindent The principal idea of the ECM is to find a set of $\amplitudal{m}$ integration points ($\amplitudal{m}\ll m$) such that
\begin{equation}
	\matr[\,k]{\amplitudal{f}}=\left[\transp{\Matr{\ROMbasis}}\cdot\int_{V}\transp{\Matr{B}}\cdot\matr{\mStress}\,\mathrm{d}V\right]^k=\left[\int_{V_0}\dfrac{1}{j}\cdot\transp{\Matr{\ROMbasis}}\cdot\transp{\Matr{B}}\cdot\matr{\mStress}\,\mathrm{d}V\right]^k\approx\sum_{g=1}^{m}w^g\matr[\,k\,g]{\amplitudal{f}}\approx\sum_{g=1}^{\amplitudal{m}}\amplitudal{w}^g\matr[\,k\,g]{\amplitudal{f}}\ ,\ \matr[\,g]{\amplitudal{f}}=\matr[]{\amplitudal{f}}(\matr[\,g]{\amplitudal{x}})
\end{equation}
\noindent for all $k=1,\dots,p$ snapshots. Preferably an integral over the reference domain should be pursued for finite deformations, so that the weights remain constant. In finding the hyper integration points among the points of the HF integration $\vect{\amplitudal{z}}=\{\vect[\,g]{\mposition}\}_{g=1}^{\amplitudal{m}}\subset\vect{z}$ and their associated positive weights $\vect{\amplitudal{w}}=\transp{[\amplitudal{w}^1,\amplitudal{w}^2,\dots,\amplitudal{w}^g]}$, an optimization problem is formulated and transferred to matrix notation to make it suitable for computer implementation, which reads
\begin{equation}\label{eqn:ECM1}
(\matr{\amplitudal{w}},\matr{\amplitudal{z}})=\arg\!\underset{\underline{\amplitudal{w}}\in\mathbb{R}^{\amplitudal{m}}_+,\,\vect{\amplitudal{z}}\subset\vect{z}}{\min}\left(\norm{\sum_{k=1}^{p}\sum_{i=1}^{\amplitudal{n}_\mathrm{t}}\sum_{g=1}^{\amplitudal{m}}\amplitudal{w}^g{\amplitudal{f}}^{k\,g}_{i}-f_{i}^{k}}^2\!\!+\norm{\sum_{g=1}^{\amplitudal{m}}\amplitudal{w}^g-V_0}^2\right)\!=\!\arg\!\underset{\underline{\amplitudal{w}}\in\mathbb{R}^{\amplitudal{m}}_+,\,\matr{\amplitudal{z}}\subset\vect{z}}{\min}\norm{\begin{bmatrix}\Matr[{\underline{\amplitudal{z}}}]{x}_f\\\underline{1}\end{bmatrix}\cdot\underline{\amplitudal{w}}-\begin{bmatrix}\underline{0}\\V_0\end{bmatrix}}\ ,
\end{equation}
\noindent with the training matrix $\Matr{x}_f$
\begin{equation}
	\Matr{x}_f=\begin{bmatrix}
		\amplitudal{f}^{1\,1}_1-\amplitudal{f}_{1}^{1}/V_0&\amplitudal{f}^{1\,2}_1-\amplitudal{f}_{1}^{1}/V_0&\cdots&\amplitudal{f}^{1\,m}_1-\amplitudal{f}_{1}^{1}/V_0\\
		\amplitudal{f}^{1\,1}_2-\amplitudal{f}_{2}^{1}/V_0&\amplitudal{f}^{1\,2}_2-\amplitudal{f}_{2}^{1}/V_0&\cdots&\amplitudal{f}^{1\,m}_2-\amplitudal{f}_{1}^{2}/V_0\\
		\cdots&\cdots&\vdots&\cdots\\
		\amplitudal{f}^{p\,1}_{\amplitudal{n}_{\mathrm{t}}}-\amplitudal{f}^{p}_{\amplitudal{n}_{\mathrm{t}}}/V_0&\amplitudal{f}^{p\,2}_{\amplitudal{n}_{\mathrm{t}}}-\amplitudal{f}^{p}_{\amplitudal{n}_{\mathrm{t}}}/V_0&\cdots&\amplitudal{f}^{p\,m}_{\amplitudal{n}_{\mathrm{t}}}-\amplitudal{f}^{p}_{\amplitudal{n_{\mathrm{t}}}}/V_0
	\end{bmatrix}%\in\mathbb{R}^{(p\cdot\amplitudal{n}_\mathrm{t})\times m}\ \ , \ \matr{b}=\matr{0}\in\mathbb{R}^{p\cdot\amplitudal{n}_{\mathrm{t}}} \ .
\end{equation} This includes the fundamental requirement of every plausible integration scheme, which states that the volume has to be integrated exactly 
\begin{equation}
	\int_{V_0}1\,\mathrm{d}V=\sum_{g=1}^{m}w^g=\sum_{g=1}^{\amplitudal{m}}\amplitudal{w}^g=V_0\ .
\end{equation}

Therein, the number $\amplitudal{n}_{\mathrm{t}}$ equaling $\amplitudal{n}_{\mathrm{t}}=\amplitudal{n}+n_{\MStress}$, reflects the total number of global force values, with $n_{\MStress}$ representing the distinct macroscopic stress components (3D:\,\,$n_{\MStress}\!=\!6$, 2D:\,\,$n_{\MStress}\!=\!3$), see equation \eqref{eqn:FEM_micro}, since in this formulation the total displacement field is included in $\Matr{\ROMbasis}$, in contrast to the implementation of Hern\'andez, where  $\Matr{\ROMbasis}$ is a basis for the fluctuation field. This ensures that the macroscopic stress field will be integrated correctly.

\citet{Hernandez2014} applied a nonnegative greedy selection algorithm to select the hyper points. In order to integrate the volume exactly and be able to efficiently carry out the selection algorithm, a SVD is applied on the training matrix $\Matr{x}_f$, instead of using it directly, as shown in the following equation
\begin{equation}\label{eqn:SVDlambda}
	\transp{\Matr{x}}_f=\Matr{\lambda}\cdot\Matr{\sigma}_{f}\cdot\transp{\Matr{v}_{f}}\ \ \text{and}\ \ \Matr{x}_f\cdot\underline{\amplitudal{w}}=\Matr{v}_{f}\cdot\transp{\Matr{\sigma}_{f}}\cdot\underbrace{\left[\transp{\Matr{\lambda}}\cdot\underline{\amplitudal{w}}\right]}_{\approx \underline{0}}\approx \underline{0}\ ,
\end{equation}  since a high degree of correlation can be observed therein.
This comes with the important side effect that $\Matr{\lambda}$ is a unitary matrix and the greedy algorithm will select an optimal result for the information given in $\Matr{\lambda}$ when it is truncated to $\amplitudal{m}-1$ rows. The overall algorithm is summarized in Box~\ref{box:ECM}.

\begin{BOX}[!h]
	\centering
	\fbox{%
		\parbox{0.8\textwidth}{
			\textbf{Inputs:} $\Matr{j}=\transp{\left[\ \Matr{\lambda}\ \,\matr{1}\ \right]}\ \ ,\ \underline{b}=\transp{\left[\ \transp{\underline{0}}\ \,V\ \right]}$\\[2ex]
			\textbf{Starting values:}\\
			- The set containing hyper integration points: $\amplitudal{\matr{z}}=\{\;\}$\\
			- The set of available points: $\matr{y}=\matr{z}$\\
			- The number of pre-chosen  (positive) weights $\amplitudal{\matr{w}}$: $\amplitudal{m}_+=0$\\
			- The residual vectors: $\matr{r}=\matr{b}$\\[2ex]
			\textbf{while} $\amplitudal{m}_+<\amplitudal{m}$ \textbf{do:}\\
			1. Compute next hyper point: $k=\arg\underset{k\in\matr{y}}{\max}\ \left\{\matr{j}_{\,k}\cdot\matr{r}/||\matr{j}_{\,k}||\right\}$ \ .\\
			2. Add $k$ to set $\amplitudal{\matr{z}}$ and remove it from $\matr{y}$ ($\amplitudal{\matr{z}}=\amplitudal{\matr{z}}\cup k$, $\matr{y}=\matr{y}\setminus k$)\ .\\
			3. Determine $\amplitudal{w}$ by unrestricted least square minimization (\texttt{numpy.linalg.lstsq}):
			\begin{center}
				$\amplitudal{\matr{w}}=\min\ \norm{\Matr[\matr{\amplitudal{z}}]{j}\cdot\underline{\amplitudal{w}}-\matr{b}}$\ .
			\end{center}
			4. If all entries of $\amplitudal{\matr{w}}$ are nonnegative, go to step 7.\\
			5. Determine $\amplitudal{\matr{w}}$ by applying nonnegative least square min. (\texttt{scipy.optimize.nnls})
			\begin{center}
				$\amplitudal{\matr{w}}=\arg\underset{\amplitudal{\matr{w}}>0}{\min}\ \norm{\Matr[\matr{\amplitudal{z}}]{j}\cdot\underline{\amplitudal{w}}-\matr{b}}$\ .
			\end{center}
			6. Set $\amplitudal{\matr{z}}=\amplitudal{\matr{z}}\setminus\amplitudal{\matr{z}}_{\,0}$ with $\amplitudal{\matr{z}}_{\,0}=\amplitudal{\matr{z}}_{\,0}\cup\{g\,|\,g\in\amplitudal{\matr{z}}\land \amplitudal{w}^g=0\}$, $\matr{\amplitudal{w}}=\matr{\amplitudal{w}}(\amplitudal{\matr{z}})$, with $\matr{z}_{\,0}$ representing a set of excluded points with zero weight.\\
			7. Compute the residual $\matr{r}=\Matr[\matr{\amplitudal{z}}]{j}\cdot\underline{\amplitudal{w}}-\matr{b}$ and set $m_+=\mathrm{card}(\amplitudal{\matr{z}})$.\\[1ex]
			\textbf{end}
		}
	}
\caption{Empirical cubature method (offline stage). Greedy algorithm for computing an optimal set of integration points $\matr{\amplitudal{z}}={1,2,\dots,\amplitudal{m}}$ and corresponding weights $\amplitudal{w}\in\mathbb{R}^{\amplitudal{m}}_+$, after \citep{Hernandez2016})}
\label{box:ECM}
\end{BOX}

\subsection{Clustered training trajectories}
The hyper ROM method relies on an \emph{offline} (training) phase to gather the matrices $\Matr{\nodal{\snapshot}}_\mdisplacement$ and $\Matr{x}_f$ including the nodal displacements and integration point forces, respectively. In this process, a set  $\mathbb{T}$ with stress $\undertilde{\Sigma}(t)$ and displacement gradient $\undertilde{H}(t)$ trajectories
\begin{equation}
    \mathbb{T}=\{\undertilde{\Sigma}_1(t),\ \undertilde{\Sigma}_2(t),\ \dots,\ \undertilde{H}_1(t),\ \undertilde{H}_2(t),\ \dots\}
\end{equation}
must be found. Finding a good training set is, however, a non-trivial task. The optimal set would include as few trajectories as possible, while the thereby generated snapshots are representative of the macroscopic problem(s), so that they can be solved within a given accuracy. In terms of the ROM method, this, in principle, leads to the smallest basis $\underline{\underline{\ROMbasis}}$, that incorporates only modes relevant for the problem(s) to be solved.

A \enquote{blind} training strategy, as, for example, applied in \citep{Logarzo2021}, may incorporate 10,000 trajectories, even in the 2D case. A more sophisticated, iterative training strategy for Neural Networks can, e.g., be found in \citep{Kalina2022}. In this approach, the macroscopic problem is simulated starting with an small initial training set. If new paths, so far not  trained, are found, the training set is enriched and after a Neural Network calibration process the problem is simulated again. This process is repeated until convergence is reached.

Since the hyper ROM method inter- and extrapolates very well in regions which were not trained, we propose a process, in which the desired structure is simulated at first using a simple surrogate material model, e.g.~only linear elastic, with the initial RVE stiffness. In force-controlled problems, the stress field computed in this manner will often be not too far off from the solution of the actual problem, and, vice versa, the estimated strain field in displacement-controlled problems. Even if large errors, which would usually be unacceptable, were to be introduced here, it would conceptually not be a problem, because they only affect the training input, which may still yield good simulation results in view of the good inter- and extrapolation properties. 

Then, the stress or strain paths at the macroscopic integration points must be clustered, to vastly reduce the amount of required training data and to minimize unnecessary recurrences. Clustering, in this context, can be understood as a type of unsupervised learning strategy, for which an overview can be found in \citep{Han2012}. We specifically make use $k$-Means Clustering, going back to \citet{Lloyd1982}, which is frequently used and widely accepted. To find the necessary number of $k$ trajectories, the Elbow Method \cite{Thorndike1953} or more sophisticated approaches may be applied. The algorithm is briefly summarized in Box \ref{box:training}.

\begin{BOX}[!h]
	\centering
	\fbox{%
		\parbox{0.8\textwidth}{
			\begin{enumerate}
			    \item Determine a suitable simple surrogate material model (e.g.~linear elastic, with initial RVE stiffness).
                \item Simulate the macroscopic structure of interest and collect the stress or strain paths at the integration points.
                \item Apply the $k$-means algorithm in suitable manner on the data to obtain $\mathbb{T}$.\label{step:kmeans}
                \item Collect the training snapshots by simulating the RVE using the determined training trajectories and evaluate the data to determine $\underline{\underline{\ROMbasis}}$,  $\amplitudal{\matr{z}}$ and $\matr{\amplitudal{w}}$.
                \item Simulate the structure(s) of interest.
			\end{enumerate}
		}
	}
\caption{Clustered training algorithm}
\label{box:training}
\end{BOX}

\section{Implementation}
The algorithmic approach described above was implemented to the commercial FE code Abaqus. The main strategy in doing so  was to use Abaqus on the macroscale and the previously established very compact, self-contained code \texttt{MonolithFE\textsuperscript{2}}, see  \citep{Lange2021_monolithic}, on the microscale. This was motivated by the fact that Abaqus is well tested and has a lot of built-in features and material routines. On the other hand, Abaqus is not reentrant---i.e.~the program cannot call itself from within---and also too inflexible to implement the hyper ROM method in a straightforward manner.

\begin{figure}[!h]
	\centering
	\includegraphics[width=\textwidth]{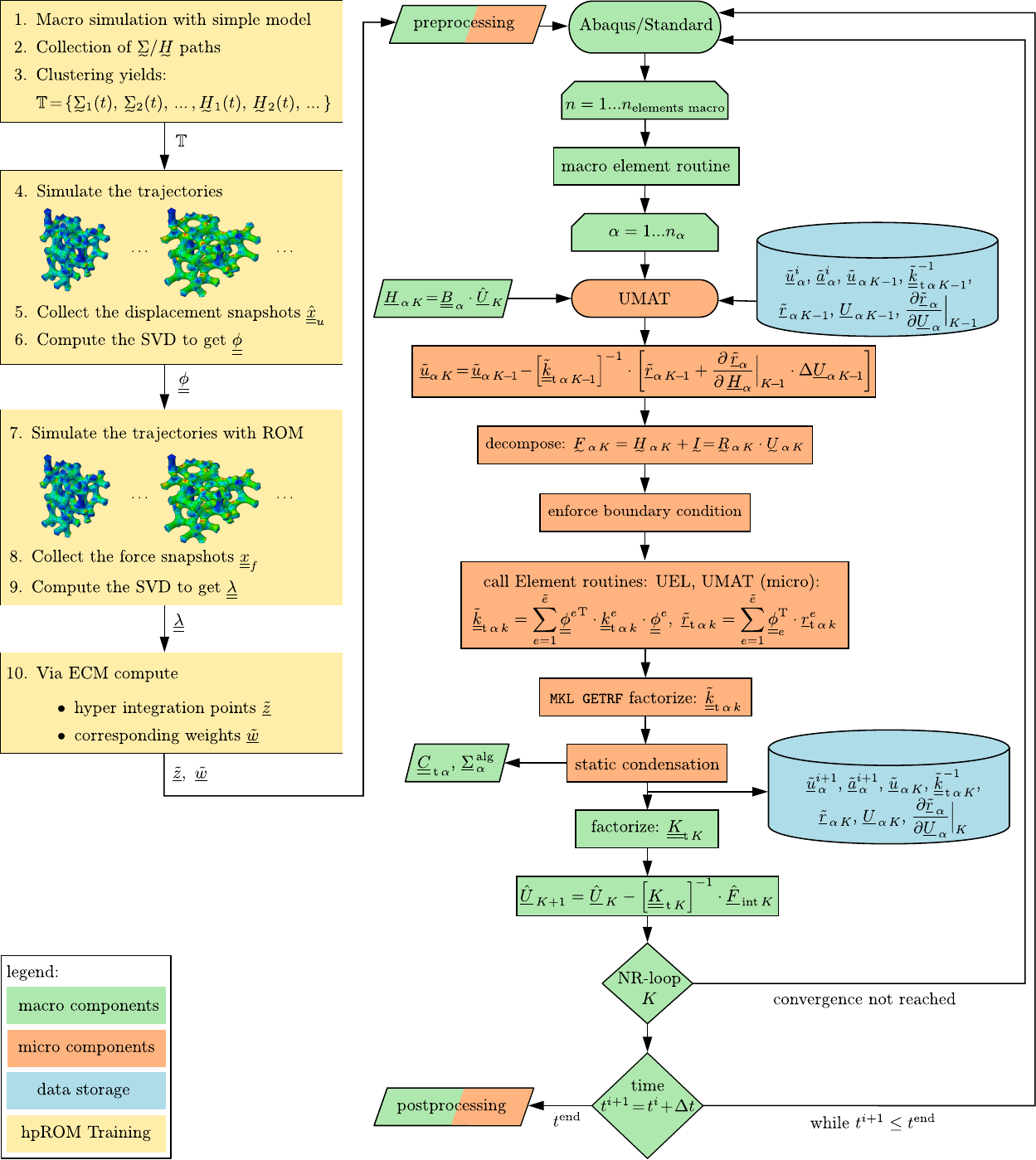}
	\caption{Flowchart with the main online and offline steps of the monolithic hyper ROM method in the MonolithFE\textsuperscript{2
 }code\label{fig:program_flow_charts}}
\end{figure}

Figure \ref{fig:program_flow_charts}  shows a flowchart that concisely illustrates the main sub-processes of the program. All parts of the \texttt{MonolithFE\textsuperscript{2}} kernel are written in object-oriented \texttt{FORTRAN 2003}. Abaqus enters the code at a specific macroscopic integration point $\alpha$ through the \texttt{UMAT} interface. Since the \texttt{UMAT} interface only supports state variables $\underline{a}_{\alpha}$ to an extent of 10,000 floating point numbers and only provides data from the last converged step, a self-written, pointer based data management had to be implemented. To make external, user-defined element and material routines compatible with Abaqus and thereby allow the testing of the respective routines independently in Abaqus, the \texttt{UEL} and \texttt{UMAT} interface is also used within \texttt{MonolithFE\textsuperscript{2}}.

The systems of equations are solved using the highly optimized \texttt{Intel \textregistered\ oneAPI Math Kernel Library} package. In HF simulations, the \texttt{PARDISO} solver is used for LU-factorizing the sparse matrix, which involves an efficient supernodal method. In ROM simulations, the \texttt{GETRF/S} solver is applied, which relies on Toledo's recursive
 LU algorithm. It is worth mentioning that the inverse operators which appear in the text and the flowcharts are only written in favor of a simple notation; inverse matrices are never actually being calculated, but rather LU factorization is utilized, together with the respective right-hand sides.

Since the FE\textsuperscript{2} method parallelizes very well, parallelization becomes a indisputable necessity regarding the high computational cost, especially for industrial applications.  The program was therefore designed to run in parallel on shared as well as on distributed memory systems.

Post-processing of FE\textsuperscript{2} computations, however, becomes a real issue, since very large data sets may occur in the HF method. Naturally, fewer data will be created in the hyper ROM method. Nevertheless, to make the quantities of interest accessible at all microscopic integration points, elaborate methods, such as Gappy Data Reconstruction \cite{Hernandez2016},  would have to be employed. From a practical perspective, and in view of the fact that often the homogenized, macroscopic values are mainly of interest, it seems more advisable to only postprocess selected RVEs, by imposing the strain history and perform a re-simulation after the actual analysis.

To generate training data, the trajectories are defined in an input file for a Driver Program, which accesses \texttt{MonolithFE\textsuperscript{2}} in parallel. In this process, the unconstrained nodal displacements $\underline{\nodal{u}}$ and the internal forces $\underline{\amplitudal{f}}$ at the integration points are written into plain-text files. Thereafter, a \texttt{Python} script conducts the hyper ROM data evaluation. The computationally intensive parts, which operate with large matrices, are written in \texttt{FORTRAN 2003} and wrapped into \texttt{Python}. Reading the raw data from the file and processing it to its final matrix representation (see $\underline{\underline{\hat{x}}}_u$ and $\underline{\underline{x}}_f$) is done, in parallel, using \texttt{openMP}. The SDV of the respective matrices is performed via the \texttt{Intel \textregistered\ oneAPI Math Kernel Library} package \texttt{GESVD}, which also can be parallelized, even on distributed-memory systems. The ECM algorithm  was written in \texttt{Python}, with aid of the \texttt{numpy} package \texttt{linalg.lstsq} and the \texttt{scipy} package \texttt{optimize.nnls} for unrestricted,  respectively nonzero-restricted, least square minimization. The results, in terms of the displacement modes $\underline{\underline{\phi}}$, hyper integration points $\underline{\amplitudal{z}}$, and corresponding weights $\underline{\amplitudal{w}}$ are written to the \texttt{MonolithFE\textsuperscript{2}} input file.

Regarding the hyper ROM method, the chosen programming style is a compromise between efficiency and keeping the necessary interventions in the classical HF code as few as possible. This preserves the modularity of the program, since otherwise a lot of code would coexist just for the sake of efficiency, although different parts, from the mathematical standpoint, would serve the same purpose. In an actual online ROM simulation, the full field of nodal displacements $\underline{\nodal{u}}$ is first recovered from the ROM amplitudes $\underline{\tilde{u}}$.

In a ROM simulation with full integration, first the full vector of internal force and corresponding stiffness matrix in the CSR format (enables sparse matrix multiplication) is assembled and afterwards projected onto the modes. In the hyper integration mode, only the elements that have at least one integration point are called and they in turn call the hyper element integration points. Only for these integration points memory is allocated to store the internal state variables $\underline{\tilde{a}}$. The full internal force vector and tangent stiffness is not build up for reasons of efficiency. Rather, the projection onto the modes follows directly with element ROM modes $\underline{\underline{\phi}}^e$, which are essentially just the rows of $\underline{\underline{\phi}}$ corresponding to respective DOFs of an element.

\section{Numerical Examples}
\subsection{General remarks}

This section will demonstrate that the monolithic hyper ROM FE\textsuperscript{2} method represents an accurate approximation of the conventional HF FE\textsuperscript{2} method, while being computationally very efficient. Accuracy is understood here in a practical engineering sense,
where deviations of around $1\%$ are considered allowable w.r.t.~the reference HF model, in terms of global (macroscale) results. Typical engineering applications usually allow even larger margins of error, since it is understood that parameters often enter the numerical analysis with substantial uncertainty. All simulations and data evaluations were carried out in parallel on a Computer with 16 threads (8 physical cores), equipped with an Intel\,\textregistered Xenon\,\textregistered Gold 6244 processor and 768 GB DDR4 RAM. Since methods presented above are modularly combinable, not all possible combinations were tested, as it was assumed that the speedup effects do not strongly influence each other.

\subsection{Example 1: Bending of a porous composite beam}

\begin{figure}[!h]
  \centering
  \subfloat[][Cantilever beam under bending, with prescribed tip displacement\label{Beam}]{\includegraphics[width=0.5\linewidth]{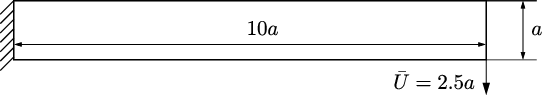}}
  \hspace{2ex}
  \subfloat[][Relative dimensions of the RVE\label{MieheRVE}]{\includegraphics[width=0.4\linewidth]{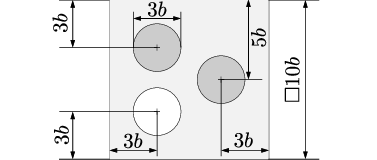}}
  \caption{Example 1: Macroscopic beam problem accounting for a micro-heterogeneous material. The RVE of the microstructure comprises a soft elastic-plastic matrix, stiff elastic inclusions and a pore\label{fig:Exp1}, cf.~\citet{miehekoch_mikro_makro_2002}. The following isotropic, linear-elastic material properties have been assumed: $E^{\mathrm{M}},\, \sigma^{\mathrm{M}}_{\mathrm{y}\,0}\!=\!0.01E^{\mathrm{M}},\,h^{\mathrm{M}}\!=\!0.016E^{\mathrm{M}}$, $E^{\mathrm{I}}\!=\!10\,E^{\mathrm{M}}$,\,$\nu^{\mathrm{M}}\!=\!\nu^{\mathrm{I}}\!=\!0.3$. Here the superscripts M and I refer to the properties of the matrix and inclusions, respectively.
  }
\end{figure}

In the first example, a two-scale simulation of the bending of a composite material beam is considered, under plain strain and small strain conditions. The beam's relative dimensions and boundary conditions are shown in Figure~\ref{fig:Exp1}(\subref{Beam}), whereas
Figure~\ref{fig:Exp1}(\subref{MieheRVE}) illustrates the assumed microstructure. 

\begin{figure}[!h]
  \centering
  \includegraphics[width=1.0\linewidth]{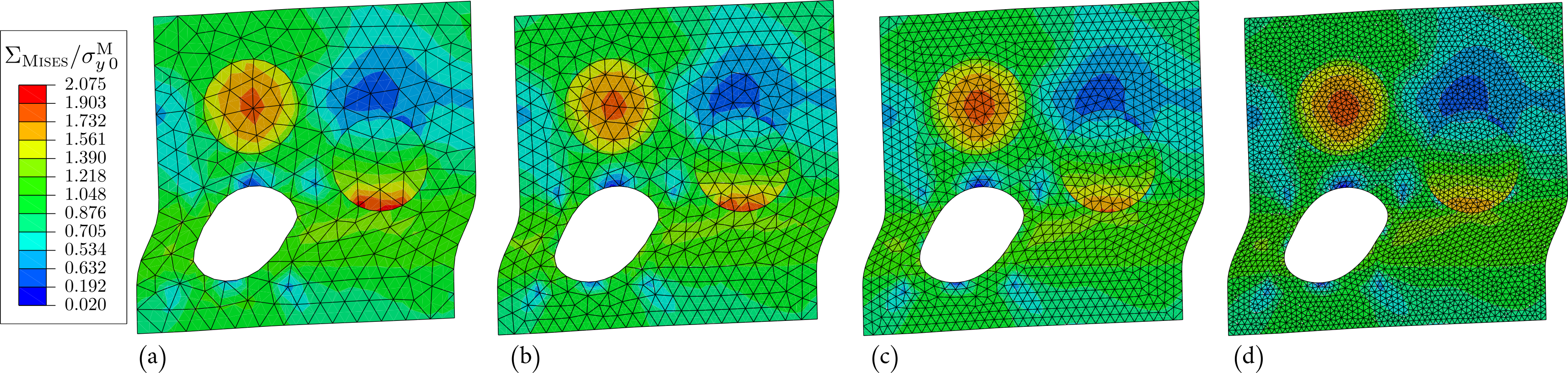}
  \caption{Meshes considered in the 2d h-refinement study and superimposed contour plots of the normalized von Mises stress fields}
  \label{fig:mesh_refinement}
\end{figure}

\begin{figure}[!h]
  \centering
  \subfloat[][Normalized stress-strain curves for meshes (a)--(d)]{\includegraphics[width=0.48\linewidth]{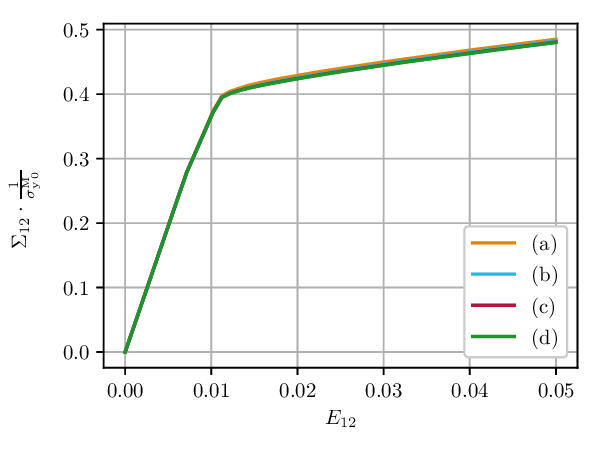}}
  \hspace{2ex}
  \subfloat[][Deviation of the stress values relative to mesh (d)]{\includegraphics[width=0.48\linewidth]{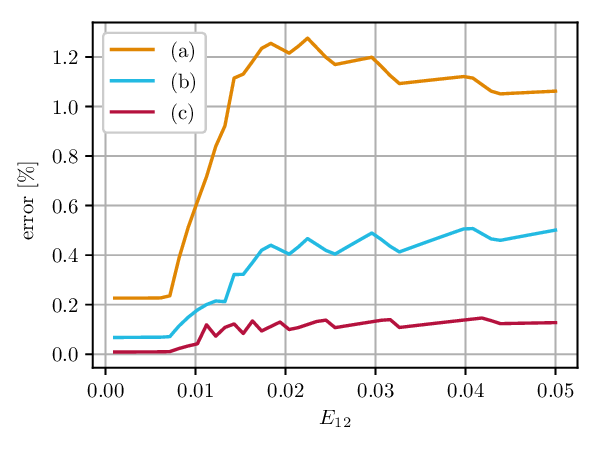}}
  \caption{Results of the mesh convergence study for Example 1 in terms of macroscopic shear stress responses}
\label{fig:mesh_refinement_curves}
\end{figure}

To apply the hyper ROM method efficiently, and for the sake of a sensible and fair comparison, a convergence study was performed first. The RVE geometry was meshed with quadratic, triangular standard elements. The specific meshes that were included in the analysis are shown in Figure \ref{fig:mesh_refinement}. A pure shear test case ($E_{11}\!=\!E_{22}\!=\!0,\,E_{12}\!=\!0.05$) was considered to determine which of the meshes have acceptable quality. The resulting macroscopic stress-strain curves are shown in Figure~\ref{fig:mesh_refinement_curves}(a). Taking mesh (d) as the reference, the relative error of the macroscopic shear stress responses was plotted in Figure \ref{fig:mesh_refinement_curves}(b). It is observed that convergence with mesh refinement is indeed achieved for the employed quadratic elements. Meshes (b) and (c) show global stress deviations relative to Mesh (d) of less than $1\%$ over the entire strain range. Taking also the convergence of the local stress field solutions into account, see contour plots in Figure~\ref{fig:mesh_refinement}, Mesh (c) was selected as the best compromise between accuracy and efficiency and therefore used in the subsequent investigations.


In the next step, training trajectories were to be chosen. In pursuing the clustered training strategy, first a rough surrogate material model had to be found. In this case, an isotropic von Mises ($J_2$) plasticity model, with power law hardening was used. It requires four material parameters, which were determined via a single, uniaxial tension simulation (with $E_{11}\!=\!0.05,\ \Sigma_{22}\!=\!\Sigma_{12}\!=\!0.0\:\mathrm{MPa}$) as $E\!=\!1.01\:E^{\mathrm{M}},\nu=0.29, \sigma_{\mathrm{y}\,0}\!=\!0.72\,\sigma^{\mathrm{M}}_{\mathrm{y}\,0}, N\!=\!0.08$. Note that the parameters were not optimized---in a classical parameter identification sense---, but rather set by approximation formulas and experience. 

Subsequently, the two-scale beam analysis was conducted with this surrogate model and the strains at all integration points extracted from the final step. Afterwards, the strain values were normalized and divided into $30$ clusters, 
using a k-means clustering algorithm. These clusters were  multiplied by the largest strain Frobenius norm of the beam model and then used as end values for the monotonic, strain-driven training simulations.

\begin{figure}[!h]
  \centering
  \subfloat[][Force-displacement curves for various ROM mode numbers $\amplitudal{n}$\label{FU_ROM_modes}]{\includegraphics[width=0.48\linewidth]{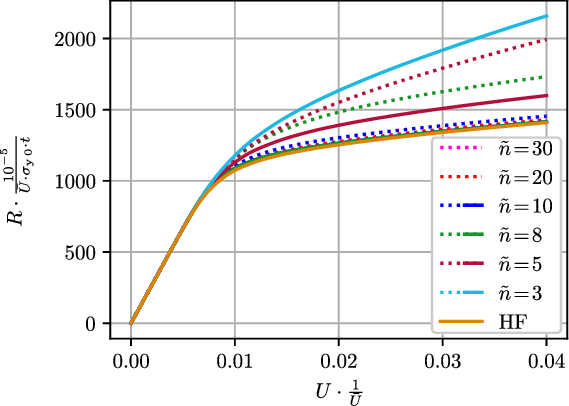}}
  \hspace{2ex}
  \subfloat[][Associated relative error w.r.t.~the HF solution \label{error_ROM_modes_wrt_HP}]{\includegraphics[width=0.48\linewidth]{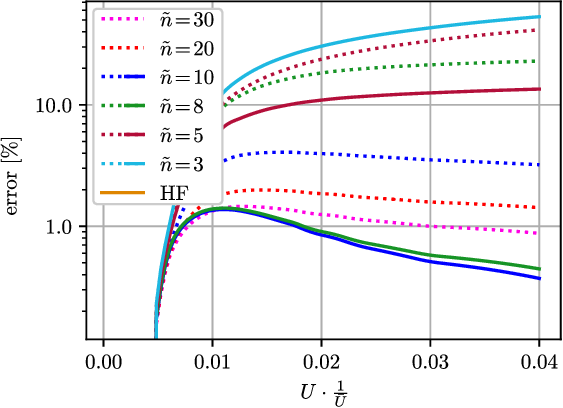}}
  \caption{ROM mode convergence study with  Clustered (solid lines) and Unspecific Training (dashed lines) \label{ROM_convergence_exp1}}
\end{figure}

To find the necessary number of ROM modes to simulate Example 1 with sufficient accuracy, the beam multiscale problem was analyzed using an increasing number of ROM modes, starting from $\tilde{n}\!=\!3$ elastic modes up to $\tilde{n}\!=\!10$ modes and the HF results for comparison. The computed global force $R$ over prescribed displacement $U$ reaction curves for the tip of the beam are displayed in Figure~\ref{ROM_convergence_exp1}\,(\subref{FU_ROM_modes}). The corresponding relative errors w.r.t.~the HF solution are shown in Figure \ref{ROM_convergence_exp1}(\subref{error_ROM_modes_wrt_HP}). It is observed that with increasing number of modes, the ROM solution quickly converges towards the HF solution. Surprisingly, $\tilde{n}\!=\!8$ modes are already capable of approximating the RVE's behavior with high accuracy, if the clustered training strategy is applied. The procedure of first simulating a problem fully and then in subsequent reduced forms to determine the minimal number of necessary ROM modes is of course not practicable for real applications, in which HF simulations might not even be feasible. It is shown simply to demonstrate the convergence behavior of the ROM method. In actual applications, experience and the testing of \enquote{representative}  trajectories not yet 
included in the training are better suited to determine the necessary number $\amplitudal{n}$. 

To show the positive effect of clustered training, an unspecific training strategy, with 125 training trajectories, was conducted in comparison. Note that including roughly four times as many training trajectories results in about four times higher simulation effort. As evident from the results in Figure~\ref{ROM_convergence_exp1}, it takes many more ROM modes to converge to the HF result in the unspecific training case. This is the expected behavior, because the modes dominant in the unspecific training set need not necessarily be the ones dominant in the actual structure. This shows empirically that the clustered training strategy not only reduces the number of training simulations, but also the number of necessary ROM modes and thus the online computational costs.

\begin{figure}[!h]
  \centering
  \subfloat[][Force-displacement curves for various  $\amplitudal{m}$\label{FU_HP}]{\includegraphics[width=0.48\linewidth]{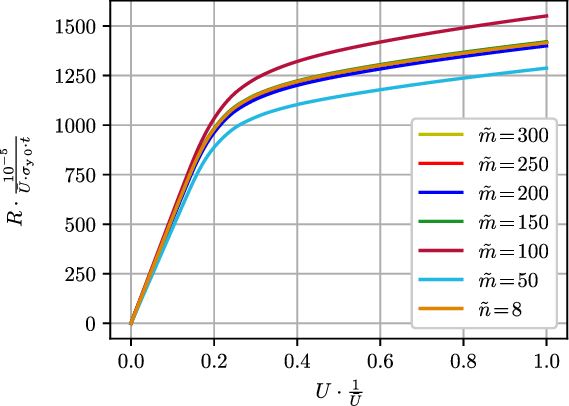}}
  \hspace{2ex}
  \subfloat[][Rel.~error w.r.t.~ROM solution for variation of $\amplitudal{m}$\label{error_HP_wrt_HF}]{\includegraphics[width=0.48\linewidth]{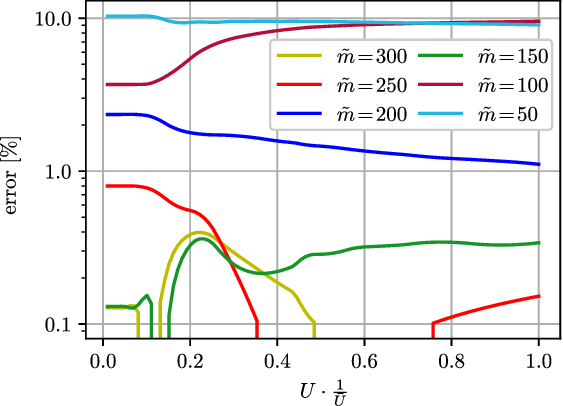}}
  \caption{Hyper integration point convergence study for a fixed number of $\tilde{n}=8$ ROM modes}
  \label{fig:HPconvergence}
\end{figure}

To find out how many hyper integration points are needed to describe Example 1 with satisfactory accuracy, another convergence study was performed. This time, the training trajectories were simulated with a restriction of the FE problem to $8$ ROM modes, while increasing the number of hyper integration points in steps of $50$, starting from the value of $\tilde{m}\!=\!50$. The corresponding global force vs.~prescribed tip displacement curves are shown in Figure~\ref{fig:HPconvergence}(\subref{FU_HP}). Deviations relative to the fully-integrated ROM results are plotted in Figure \ref{fig:HPconvergence}(\subref{error_HP_wrt_HF}). The results show a convergence towards the reference solution with an increasing number of integration points, albeit not in a monotonous fashion. The hyper 8-mode ROM model with $\tilde{m}\!=\!250$ is associated with a relative error below $1\%$ (globally), and was therefore accepted as the final hyper ROM model. 

\begin{figure}[!h]
  \centering
  \subfloat[][von Mises stress fields at the macro- and, for one selected macro integration point, also the microscale,  at $U=\bar{U}$ (HF solution)\label{Mises_stress_example1}]{\includegraphics[width=0.48\linewidth]{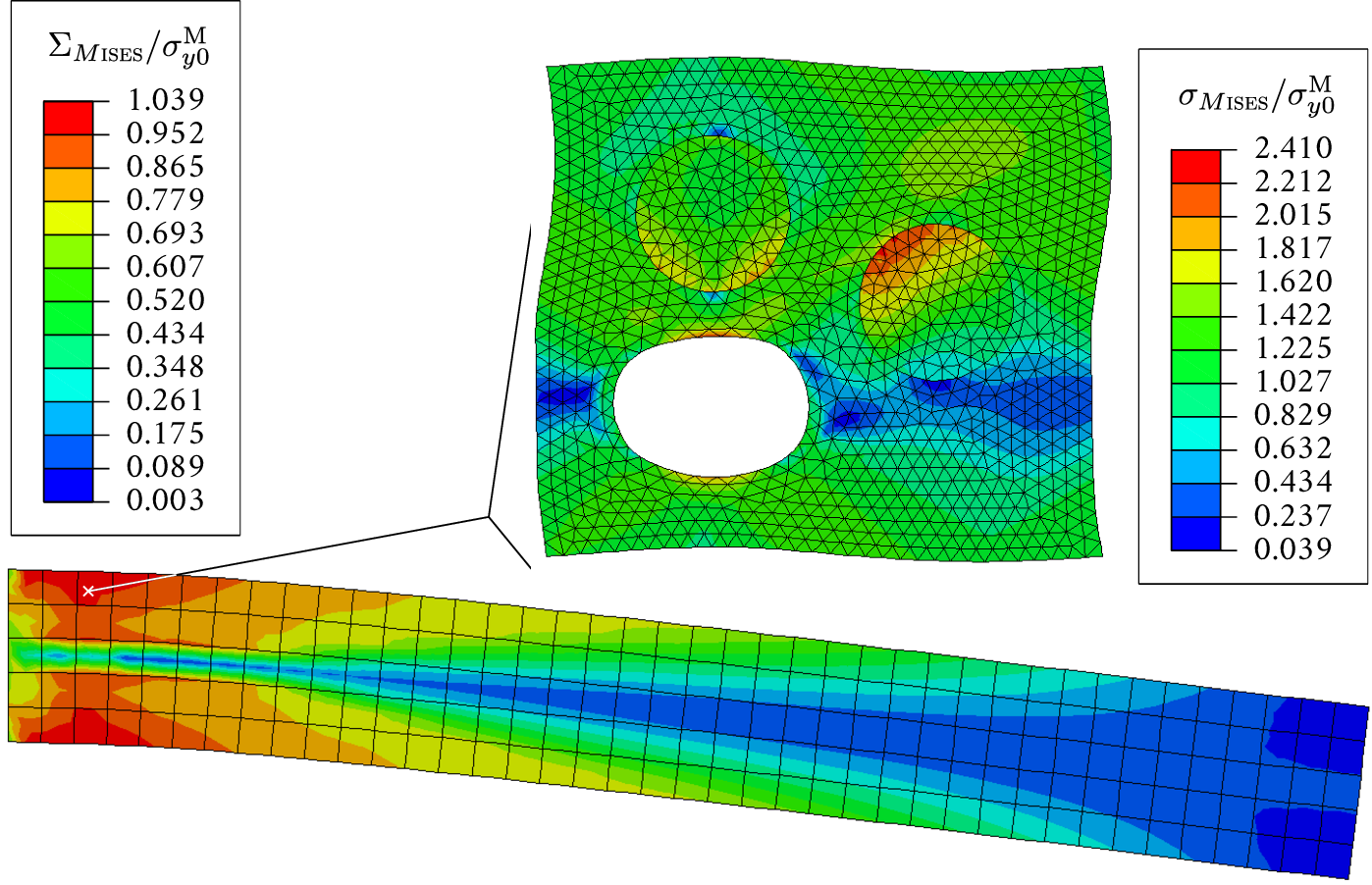}}
  \hspace{2ex}
  \subfloat[][Computed force-displacement curves with elastic recovery, for the HF, ROM and hyper ROM FE\textsuperscript{2} methods \label{UF_w_elastic_recovery}]{\includegraphics[width=0.48\linewidth]{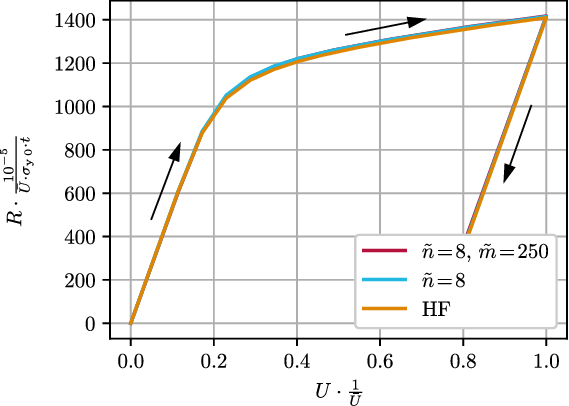}}
  \caption{Final two-scale simulation result for beam bending of a composite}
\label{fig:finalresultexample1}
\end{figure}

With the number of ROM modes and hyper integration points fixed, the multiscale analysis of the composite beam problem was then performed for the full loading sequence. This  consisted of ramping up the prescribed displacement to the maximum value of $\bar{U}\!=\!2.5\,a$ and then unloading the beam until the global tip reaction force $R$ returned back to zero. Figure \ref{fig:finalresultexample1}(\subref{Mises_stress_example1}) depicts the von Mises stress field throughout the beam at peak load, while the insert shows the microscopic equivalent stress distribution for a selected macro integration point. The global reaction force over  prescribed tip displacement curves as predicted by the different algorithms are presented in Figure~\ref{fig:finalresultexample1}(\subref{UF_w_elastic_recovery}). It turns out that the ROM and hyper ROM FE\textsuperscript{2} simulations predict the correct irreversible behavior, despite the fact that only monotonic training cases were considered.

\begin{figure}[ht!]
    \centering
  \includegraphics[width=0.7\linewidth]{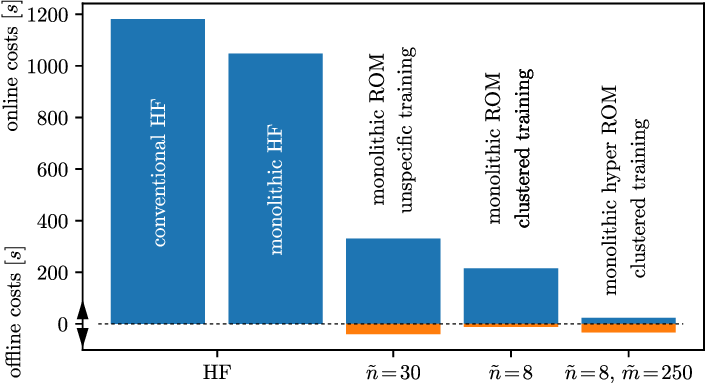}
   \caption{Computational effort associated with the different numerical schemes for Example 1}
   \label{fig:simulation_time_example_1}
\end{figure}

\begin{table}[ht!]
\centering
\begin{tabular}{|p{0.14\linewidth}|p{0.16\linewidth}|p{0.16\linewidth}|p{0.1\linewidth}|p{0.1\linewidth}|p{0.12\linewidth}|}
  \hline
  & conventional HF & monolithic HF &  \multicolumn{2}{|p{0.2\linewidth}|}{ROM} & hyper ROM \\\hline
  training& \multicolumn{2}{|c|}{-}&unspecific&\multicolumn{2}{|c|}{clustered}\\\hline
  DOFs & \multicolumn{2}{|c|}{16\,385\,400\,(100\%)} & \multicolumn{1}{|c|}{54\,000\,(<1\%)}&\multicolumn{2}{|c|}{14\,400\,(<1\%)}\\\hline
  IPs & \multicolumn{4}{c|}{12\,209\,400\,(100\%)} &\multicolumn{1}{c|}{450\,000\,(4\%)}\\\hline
  online time & \multicolumn{1}{c|}{1\,181\,s (100\%)}  & \multicolumn{1}{c|}{1\,048\,s (89\%)}& \multicolumn{1}{c}{330\,s (28\%)}  & \multicolumn{1}{c}{215\,s (18\%)} & \multicolumn{1}{|c|}{23\,s (1.9\%)}\\\hline
  offline ROM & \multicolumn{2}{|c|}{-} &\multicolumn{1}{|c|}{40\,s (3.4\%)}& \multicolumn{2}{|c|}{12\,s (1.0\%)} \\\hline
  offline hyper& \multicolumn{4}{c|}{-} &\multicolumn{1}{c|}{21\,s (1.8\%)}\\\hline
  \end{tabular}
\caption{Simulation effort in Example 1: ROM offline time includes simulation of the beam with a surrogate model, RVE-FE training  simulations and SVD on $\underline{\underline{\hat{x}}}_u$, hyper offline time includes ROM-RVE-FE training  simulations, SVD on $\underline{\underline{x}}_f$ and ECM selection\label{table:simulation_times_composite}}
\end{table}

Since the main objective of the hyper ROM method is to reduce the computational effort in terms of \emph{total} required simulation time, a more detailed analysis of its various on- and offline contributions is of interest here. The results are illustrated in Figure \ref{fig:simulation_time_example_1}, in which the online time is shown in blue and the offline time in orange. More details are given in Table \ref{table:simulation_times_composite}. 

It can be seen that the monolithic hyper ROM FE\textsuperscript{2} vastly reduced the online simulation time of the FE\textsuperscript{2} problem in Example 1 to less than $2\%$ of the high fidelity analysis and to still only about $5\%$, if the total computational effort is considered. One further observes that unspecific training requires more than three times the training effort of the clustered training. This is due to the higher number of modes required to match the accuracy requirements, resulting in about $50\%$ more online time.

\subsection{Example 2: Ceramic flow-through filter}
\begin{figure}[!h]
  \centering
  \subfloat[][Relative filter dimensions: $R/H\!=\!1.5$, material parameters: $E$, $\nu\!=\!0.14$, $A\!=\!22.09\,E$,  $n\!=\!1.06$, $m\!=\!-0.56$, cf.~\citep{Lange2022_filter})\label{fig:Filter_sketch}]{\includegraphics[width=0.70\linewidth]{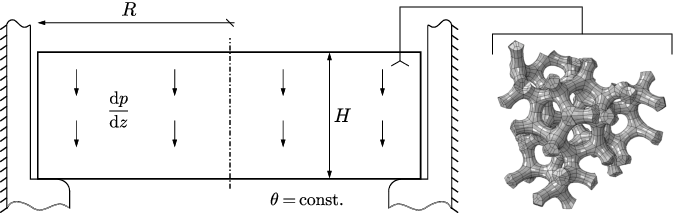}}
  \hspace{2ex}
  \subfloat[][Loading curve \label{fig:Filter_loading}]{\includegraphics[width=0.27\linewidth]{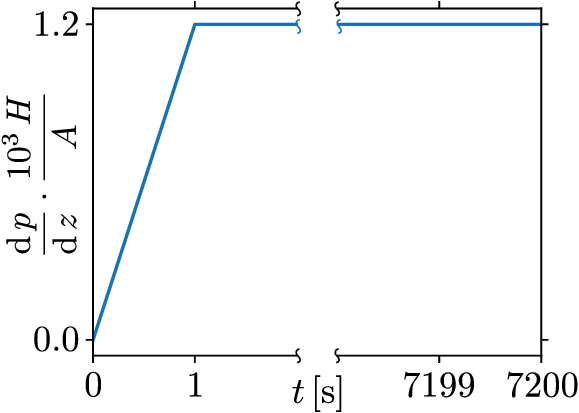}}
  \caption{Example 2: Ceramic flow-through foam filter at high temperatures, examination of creep behavior for 2h of metal flowing\label{fig:problem_creep_filter}}
\end{figure}

For the second example, a flow-through ceramic filter is considered as depicted in Figure \ref{fig:problem_creep_filter}(\subref{fig:Filter_sketch}). Such open-cell foam structures are used in continuous casting processes to filter out unwanted constituents from the metal melt. The through-flowing melt acts as a volumetric body force, from a mechanical point of view. For more details see \citet{Lange2022_filter}. On the macroscale, an axisymmetric model with linear, fully-integrated Abaqus elements (CAX4) at large deformations are used. The contact between filter and coating is described by Coulomb friction, with a coefficient of $\mu\!=\!0.3$. The mechanical loading---in terms of a homogeneously distributed local pressure drop  over time---is  linearly increased within the first second and then held constant for two hours,  as visualized in Figure \ref{fig:problem_creep_filter}(\subref{fig:Filter_loading}).
\\
The RVE (with relative density $\varrho_\mathrm{rel}\!=\!0.20$ and strut form factor $k\!=\!0.75,\, \varrho_\mathrm{rel}\!=\!0.20$) and its mesh were taken from  \citep{Lange2022_filter}. The material behavior is represented through a hypo-elastic creep model in strain hardening form. Here, the Jaumann rate of the Cauchy stress tensor is specifically given by
\begin{equation}
    \stackrel{\nabla}{\undertilde{\sigma}}=3K\cdot \mathrm{sph}(\undertilde{d})+2\mu\left[\,\mathrm{dev}(\undertilde{d})-\dfrac{3}{2}\dfrac{\mathrm{dev}(\undertilde{\sigma})}{\sigma_\mathrm{eq}}\cdot\dot{\varepsilon}^\mathrm{cr}_\mathrm{eq}\,\right],\ \ \dot{\varepsilon}^\mathrm{cr}_\mathrm{eq}=\left[\dfrac{\sigma_\mathrm{eq}}{A}\right]^{\frac{n}{m+1}}\left[\,[m+1]\,\varepsilon^\mathrm{cr}_\mathrm{eq}\,\right]^{\frac{m}{m+1}},\ \ \sigma_\mathrm{eq}=\sqrt{3\,J_2}\ .
    \label{eqn:creep_law}
\end{equation}
The equations are discretized in time by means of the implicit Euler scheme and solved using the Newton-Raphson algorithm. Globally, the Updated Lagrangian method is applied by using the configuration of the last time step as reference and evaluating the rate of deformation tensor $\undertilde{d}$ through the midpoint rule after \citet{HughesWinget1980}.
\begin{figure}[!h]
  \centering
  \subfloat[][Macroscopic, circumferential stress field at $t\!=\!7200\,\mathrm{s}$ and microscopic von Mises stress distribution in one selected RVE]{\includegraphics[width=0.55\linewidth]{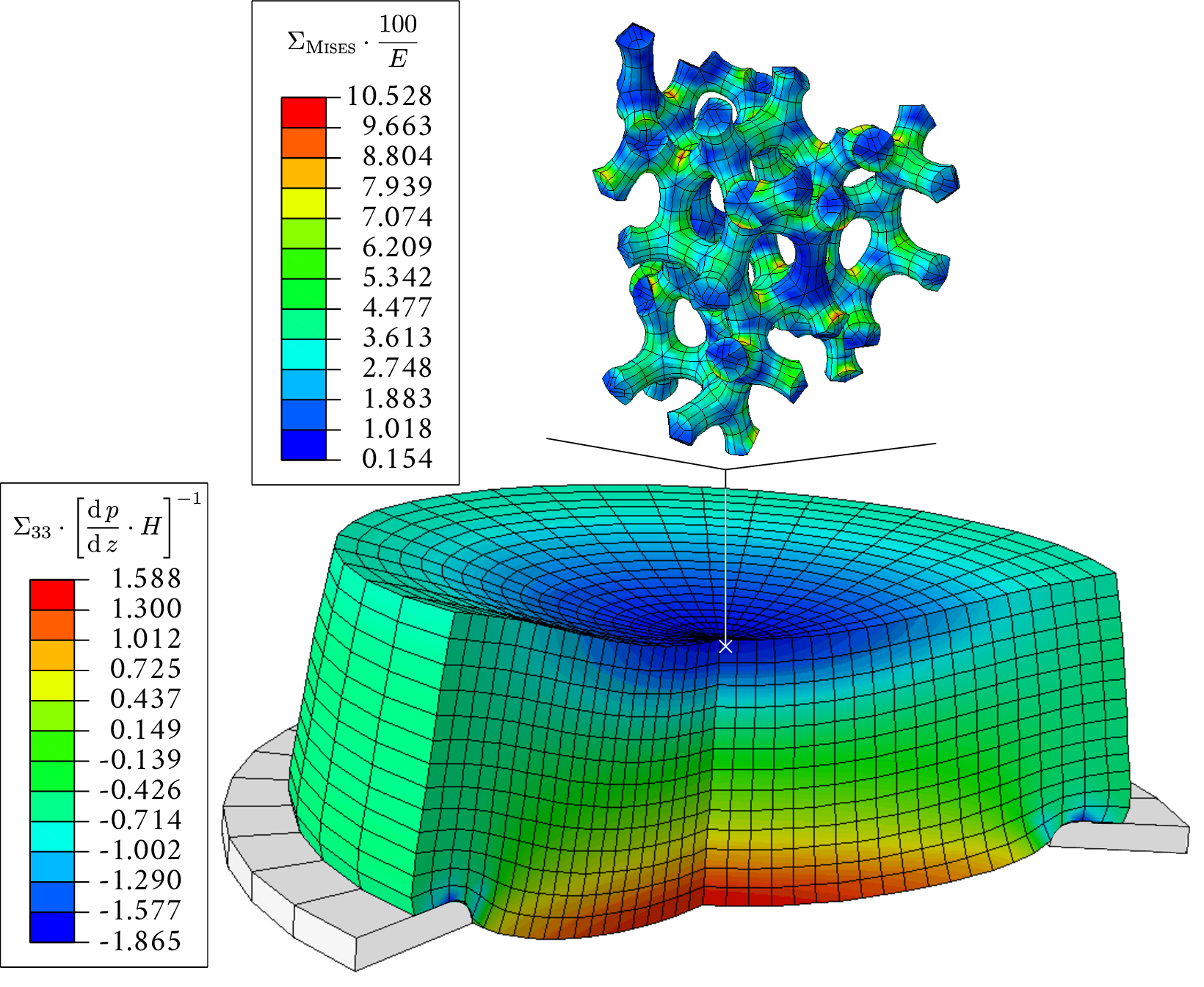}}
  \hspace{2ex}
  \subfloat[][Displacement evolution (at bottom center of the filter) over time for the HF, ROM and hyper ROM FE\textsuperscript{2} methods \label{tu_foamfilter}]{\includegraphics[width=0.42\linewidth]{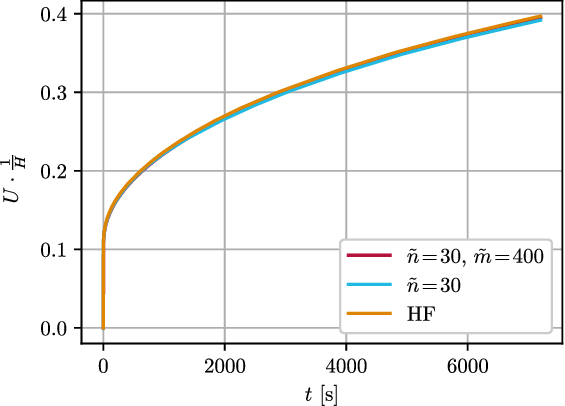}}
  \caption{Final simulation result for the (high temperature) creep behavior of a foam-like ceramic flow-through filter}
\label{fig:simulation_results_filter}
\end{figure}

At first, a simulation of the macroscopic problem with \enquote{effective} properties, also based on the constitutive law summarized in \eqref{eqn:creep_law} is conducted. Therein, the symmetry of initial effective elastic stiffness tensor is assumed to be orthotropic in nature and statically condensed out  of the RVE-FE model. The effective material parameter $A$ is fitted using a simple virtual tension test. The  stress states at the macroscopic integration points at time $t=7200\,\mathrm{s}$ is ordered into 30 clusters. Afterwards, the training procedure is set up using the stress clusters with the respective loading paths following the same pattern as the actual global loading profile shown in Figure \ref{fig:problem_creep_filter}(\subref{fig:Filter_loading}). From the resulting training data, $\tilde{n}\!=\!30$ modes are extracted. Then, the training is repeated using the ROM method and $\tilde{m}\!=\!500$ hyper integration points are extracted via the ECM algorithm.

\begin{figure}[!h]
\centering
  \includegraphics[width=0.5\linewidth]{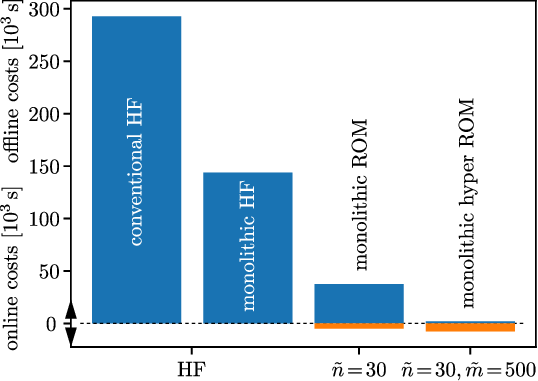}
   \caption{Computational effort associated with the different numerical schemes for Example 2}
   \label{fig:simulation_time_example_2}
\end{figure}

Figure \ref{fig:simulation_results_filter} illustrates the simulation results, where the macroscopic circumferential stress field is depicted on the left, with an insert showing the 
von Mises micro-stress field for a selected RVE. The right picture shows the macroscopic $U$-$t$-curve for all employed simulation methods. All three response curves show a good agreement---also in terms of the local measures, not pictured here. Again, the on- and offline simulation efforts are compared, see Figure~\ref{fig:simulation_time_example_2} and Table~\ref{table:simulation_times_foam}. The simulation effort in terms of  online time amounts to $13\%$ for the ROM and to below $1\%$ in the hyper ROM case. When the offline effort are additionally considered, the total computational costs are $15\%$ and $3\%$ for the ROM and hyper ROM methods, respectively, as again compared to the conventional HF model. An unspecific training was not conducted in this case, since it was clear that with the same training density as in Example 1 and four independent strain components (axisymmetry), now 625 training trajectories would have been needed, resulting in about 21 times the effort of the clustered training.

\begin{table}[ht!]
    \centering
\begin{tabular}{|p{0.12\linewidth}|p{0.16\linewidth}|p{0.16\linewidth}|p{0.14\linewidth}|p{0.14\linewidth}|}
  \hline
  & conventional HF & monolithic HF & ROM & hyper ROM \\\hline
  DOFs & \multicolumn{2}{|c|}{51\,858\,560\,(100\%)} & \multicolumn{2}{|c|}{31\,200(<1\%)}\\\hline
  IPs & \multicolumn{3}{c|}{93\,000\,960\,(100\%)} &\multicolumn{1}{|c|}{520\,000\,(<1\%)}\\\hline
  online time & \multicolumn{1}{c|}{292\,819\,s (100\%)} & \multicolumn{1}{c|}{143\,837\,s (49\%)} & \multicolumn{1}{c}{37\,551\,s (13\%)} & \multicolumn{1}{|c|}{2\,077\,s (0.7\%)}\\\hline
  offline ROM & \multicolumn{2}{|c|}{-} & \multicolumn{2}{|c|}{4\,959\,s (1.7\%)} \\\hline
  offline hyper & \multicolumn{3}{c|}{-}&\multicolumn{1}{c|}{2\,550\,s (0.9\%)}\\\hline
  \end{tabular}
\caption{Simulation effort in Example 2: ROM offline time includes simulation of the filter with a surrogate model, RVE-FE training  simulations and SVD on $\underline{\underline{\hat{x}}}_u$, hyper offline time includes ROM-RVE-FE training  simulations, SVD on $\underline{\underline{x}}_f$ and ECM selection
}\label{table:simulation_times_foam}
\end{table}

\clearpage
\subsection{Example 3: Woven composite plate with a hole}
\begin{figure}[!h]
  \centering
  \subfloat[][Plate with a hole and rotated basis system at the microscale\label{fig:Plate_w_hole_and_RVE}]{\includegraphics[width=0.38\linewidth]{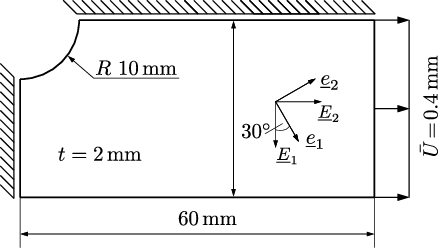}}
  \hspace{1ex}
  \subfloat[][Woven composite RVE, meshed with $109\,287$ linear tetrahedrons\label{fig:Plate_w_hole_and_RVE}]{\includegraphics[width=0.32\linewidth]{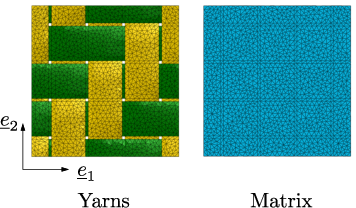}}
  \hspace{1ex}
  \subfloat[][Loading curve \label{fig:plate_w_hole_loading}]{\includegraphics[width=0.27\linewidth]{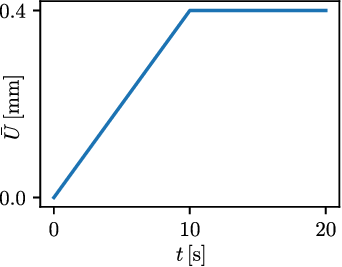}}
  \caption{Example 3: Plate with hole of a woven composite with glass fibres and thermoplastic matrix\label{fig:problem_plate_w_hole}}
\end{figure}

As a third example a tension--relaxation test of a plate with hole (quarter model) as shown in Figure \ref{fig:problem_plate_w_hole} is considered, whose microstructure is a woven glass fabric with polyamide (PA) 66 matrix, where the weft orientation is rotated by $30^\circ$ w.r.t.~the plate edge about the thickness direction. The macro problem, the RVE and the material formulation and parameters are taken from \citet{Tikarrouchine2021}. Woven composites are highly relevant engineering materials, since they have a high strength to density ratio and their behavior can be tailored to an application in the design process. If loaded into the direction of the yarns (as usually intended), the material response is anisotropic but approximately elastic, since the yarns carry most of the load. If in plane shear occurs, the matrix will also carry a considerable share of the load and the overall material behavior will become quite complex and irreversible, since PA 66 shows itself damage, viscoelasticity and viscoplasticity effects. Further complexity is introduced, because cracks and debonding effects occur in the matrix-yarn phase. This leads to very complex material formulations, which require together with a 3D--RVE a tremendous amount of computational resources, which underlines the necessity of appropriate reduction techniques.

The matrix formulation accounts for viscoelasticity, viscoplasticity and damage in the small strain setting. Details can be found in \citet{Praud2017}. As shown in \eqref{eqn:Matrix_stress}, the strain is additively split into an elastic, viscoelastic and viscoplastic share. The elastic energy release rate is denoted by $Y$ and the hardening $R$ is given in form of a power law.
\begin{equation}
    \undertilde{\sigma}=[1-d]\,\undertilde{\undertilde{C}}_\mathrm{e}:\undertilde{\varepsilon}_\mathrm{e},\ \ \undertilde{\varepsilon}_\mathrm{e}=\undertilde{\varepsilon} - \undertilde{\varepsilon}_\mathrm{p} - \sum_{i=1}^{N}\undertilde{\varepsilon}_{\mathrm{v}i},\ \ \left\{\undertilde{\sigma}_{\mathrm{v}i}=\undertilde{\undertilde{C}}_{\mathrm{v}i}:\undertilde{\varepsilon}_{\mathrm{v}i}\right\}_{i=1}^{N},\ \ Y=\dfrac{1}{2}\left[\dfrac{\undertilde{\sigma}:\undertilde{\varepsilon}_\mathrm{e}}{1-d}+\sum_{i=1}^{N}\undertilde{\sigma}_{\mathrm{v}i}:\undertilde{\varepsilon}_{\mathrm{v}i}\right],\ \ R=K\,r^n
    \label{eqn:Matrix_stress}
\end{equation}
The evolution of the internal state variables is shown in equation \eqref{eqn:Matrix_STATEV}. The plastic multiplier $r$ begins to evolve only if the yield function $f$ is positive. The evolution of the damage variable $d$ and the plastic strain $\undertilde{\varepsilon}_\mathrm{p}$ is coupled with the evolution of $r$. The viscoelastic branches are always active.
\begin{equation}
    f=\dfrac{\sqrt{3J_2}}{1-d}-R-R_0,\ \ \dot{r}=\dfrac{\langle f\rangle_+}{\mathrm{H}^{1/\mathrm{m}}},\ \ \dot{d}=\left[\dfrac{Y}{\mathrm{S}}\right]^\beta\dfrac{\dot{r}}{1-d},\ \ \undertilde{\dot{\varepsilon}}_\mathrm{p}=\dfrac{3\,\mathrm{dev}(\undertilde{\sigma})\,\dot{r}}{2[1-d]\sqrt{3J_2}},\ \ \left\{\dot{\undertilde{\varepsilon}}_\mathrm{vi}=\underline{\underline{V}}_{\mathrm{v}i}^{-1}:\left[\dfrac{\undertilde{\sigma}}{1-d}-\undertilde{\sigma}_{\mathrm{v}i}\right]\right\}_{i=1}^{N}
    \label{eqn:Matrix_STATEV}
\end{equation}
To describe the yarns, which consists of unidirectional fibres with matrix material in between, it is assumed that viscous effects can be neglected, since the matrix share amounts to only $\approx15\%$. Details can be found in \citep{Praud2017b}. The material formulation describes the formation of micro cracks with volumetric share $\gamma_c$ and related inelastic strains $\undertilde{\varepsilon}_\mathrm{ir}$. The Mori-Tanaka method is used, to describe the fictitious material \enquote{homogeneous \enquote{matrix} with cracks}. The crack shape is assumed to be known and shall be ellipsoidal, which gives a constant Eshelby tensor $\undertilde{\undertilde{S}}$, from which with the interaction tensor $\undertilde{\undertilde{T}}_c$ the actual strain concentration tensor $\undertilde{\undertilde{A}}_0$ can be computed, giving the strain in the fictions \enquote{matrix}, from which one can compute the overall stress $\undertilde{\sigma}$ as shown in \eqref{eqn:Yarn_stress}. 
\begin{equation}
    \undertilde{\sigma}=[1-\gamma_c]\,\undertilde{\undertilde{C}}_0:\undertilde{\undertilde{A}}_0:\undertilde{\varepsilon}_\mathrm{e},\ \ \undertilde{\varepsilon}_\mathrm{e}=\undertilde{\varepsilon} - \undertilde{\varepsilon}_\mathrm{ir},\ \ \undertilde{\undertilde{A}}_0=\left[\undertilde{\undertilde{I}}+\gamma_c\,\left[\undertilde{\undertilde{T}}_c-\undertilde{\undertilde{I}}\right]\right]^{-1},\ \ \undertilde{\undertilde{T}}_c=\left[\undertilde{\undertilde{I}}-\undertilde{\undertilde{S}}\right]^{-1}
    \label{eqn:Yarn_stress}
\end{equation}
The microcracks evolve when the value of a stress dependent indicative function $H_c$ is larger than $1$ and its ever encountered value at a material point as shown in equation \eqref{eqn:Yarn_STATEV}. The evolution of the inelastic strains $\undertilde{\dot{\varepsilon}}_\mathrm{ir}$ is coupled to this evolution and controlled by another indicative function $H_s$. The function $H_c$ accounts for the fact that the cracks only evolve, when the yarns are sheared or pulled against the longitudinal axis.
\begin{equation}
    \gamma_c=\gamma_c^\infty\left[1-\mathrm{exp}({-\left[\dfrac{\langle\mathrm{sup}(H_c)-1\rangle_+}{S}\right]^\beta})\right],\ \ H_c=\sqrt{\undertilde{\sigma}:\undertilde{\undertilde{H}}:\undertilde{\sigma}},\ \ \undertilde{\dot{\varepsilon}}_\mathrm{ir}=\dfrac{\undertilde{\undertilde{F}}:\undertilde{\sigma}}{H_s}\dot{\gamma}_c,\ \ H_s=\sqrt{\undertilde{\sigma}:\undertilde{\undertilde{F}}:\undertilde{\sigma}}
    \label{eqn:Yarn_STATEV}
\end{equation}
Both material formulations were discretized in time by the Euler implicit method and the resulting equations are solved using a full Newton-Raphson scheme, in contrast to convex cutting plane, general return approach presented in \citep{Praud2017,Praud2017b}.
\begin{figure}[!h]
  \centering
  \subfloat[][False color plot of the macroscopic M\textsc{ises} stress, the damage $d$ distribution in the matrix phase and microcrack density $\gamma_c$ in the yarn phase for a selected RVE\label{fig:exp3_falsecolor}]{\includegraphics[width=0.55\linewidth]{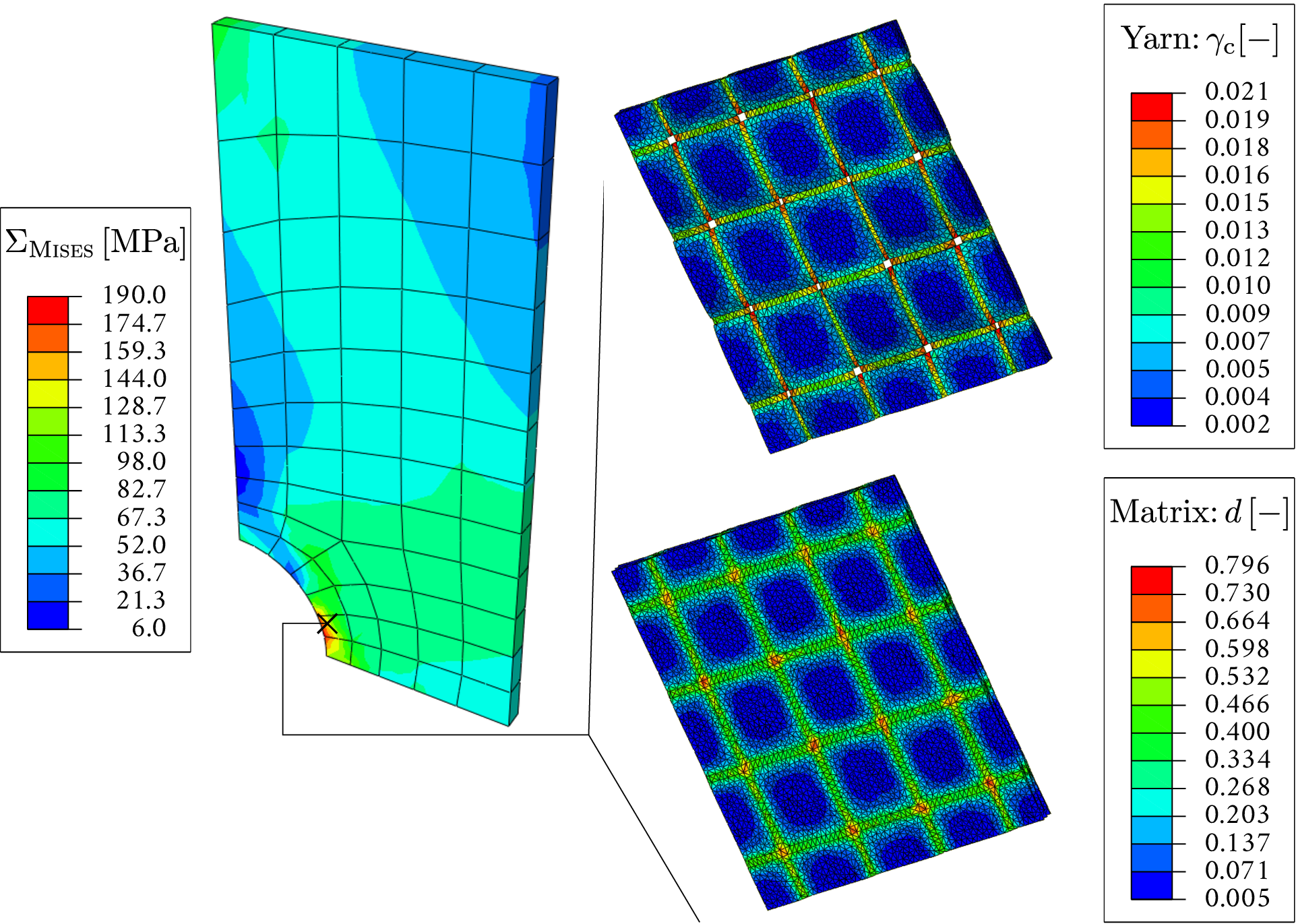}}
  \hspace{2ex}
  \subfloat[][Time force curves for the HF, ROM and
hyper ROM FE\textsuperscript{2} method\label{fig:F_t_curves_exp3}]{\includegraphics[width=0.42\linewidth]{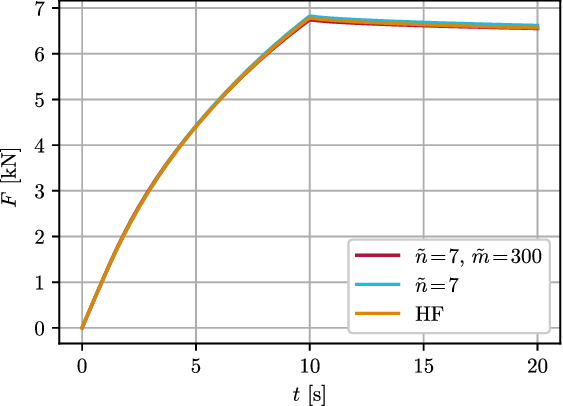}}
  \caption{Final simulation result for the tension--relaxation--test of a plate with hole made up of a woven composite with polyamide 66 matrix and glass fibres}
  \label{fig:exp3_results}
\end{figure}

To get suitable training directions for the ROM training, the initial elastic stiffness of the RVE is extracted and used as an anisotropic elastic surrogate material. From the simulation of the macro structure, meshed with quadratic, reduced integrated Abaqus C3D20R elements, 16 strain clusters were considered sufficient and were used as monotonic training paths throughout a time span of $10\,\mathrm{s}$. From the displacement snapshot matrix $\tilde{n}\!=\!7$ modes were extracted, whereby the elbow method was used. Using the ECM method $\tilde{m}\!=\!300$ integration points and corresponding weights were determined, whereby the number was chosen by experience.

\begin{figure}[!h]
\centering
  \includegraphics[width=0.5\linewidth]{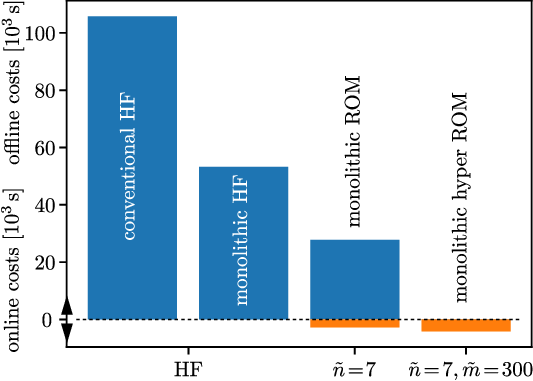}
   \caption{Computational effort associated with the different numerical schemes for Example 3\label{fig:barplot_exp3}}
\end{figure}

\begin{table}[ht!]
    \centering
\begin{tabular}{|p{0.12\linewidth}|p{0.16\linewidth}|p{0.16\linewidth}|p{0.14\linewidth}|p{0.14\linewidth}|}
  \hline
  & conventional HF & monolithic HF & ROM & hyper ROM \\\hline
  DOFs & \multicolumn{2}{|c|}{28\,124\,672\,(100\%)} & \multicolumn{2}{|c|}{3\,584(<1\%)}\\\hline
  IPs & \multicolumn{3}{c|}{55\,950\,336\,(100\%)} &\multicolumn{1}{|c|}{153\,600\,(<1\%)}\\\hline
  online time & \multicolumn{1}{c|}{105\,782\,s (100\%)} & \multicolumn{1}{c|}{53\,213\,s (50\%)} & \multicolumn{1}{c}{27\,722\,s (26\%)} & \multicolumn{1}{|c|}{110\,s (0.1\%)}\\\hline
  offline ROM & \multicolumn{2}{|c|}{-} & \multicolumn{2}{|c|}{2\,792\,s (2.6\%)} \\\hline
  offline hyper & \multicolumn{3}{c|}{-}&\multicolumn{1}{c|}{1\,361\,s (1.3\%)}\\\hline
  \end{tabular}
\caption{
Simulation effort in Example 3: ROM offline time includes simulation of the plate with a hole with a surrogate model, RVE-FE training  simulations and SVD on $\underline{\underline{\hat{x}}}_u$, hyper offline time includes ROM-RVE-FE training  simulations, SVD on $\underline{\underline{x}}_f$ and ECM selection\label{table:simulation_times_exp3}} 
\end{table}

Figure \ref{fig:exp3_results}(\subref{fig:exp3_falsecolor}) shows the distribution of the macroscopic M\textsc{ises} stress at $t\!=\!10\,\mathrm{s}$, as well as the damage $d$ in the matrix and the micro crack density $\gamma_c$ in the yarns. From the force--displacement curves in Figure \ref{fig:exp3_results}(\subref{fig:F_t_curves_exp3}) it can be seen that the ROM and hyper ROM method are in good agreement with the full solution. The barplot of Figure \ref{fig:barplot_exp3} shows that the hyper ROM method saves a lot of computational effort compared to the conventional method, whereby a speedup of 25 is achieved including the training effort and 961 if only the online time is considered, cf.~Table \ref{table:simulation_times_exp3}. The high speedup stemming from the hyper integration comes from the fact that in this example the material routine calls are extremely costly and a reduction thereof has a huge impact.

%\clearpage
\section{Summary and Conclusions}
The conventional formulation of the FE\textsuperscript{2} method represents a very flexible tool for multiscale simulations. The enormous costs are a hindrance for its application in the engineering environment. The hyper ROM method seems to be a promising approximation technique, because it still solves the original PDE, while the results lie in reduced spaces, since a lot of correlation can be observed in the HF model. The main advantage therefore is that a lot of the flexibility of the conventional HF FE\textsuperscript{2} method is preserved, while the computational costs are reduced drastically. In comparison to common hybrid NN approaches, the online costs will be higher, however much less training effort is needed and the method promises more flexibility.

A formulation based on the work of \citet{Hernandez2014} has been presented and combined with a monolithic solution scheme. The necessary amount of training data was vastly reduced by clustering loading paths from the structure of interest and eliminating rigid body rotations from the training set using a polar decomposition in the micro-macro transition rule. Numerical examples prove the applicability for practical problems and a speedup factor of almost 1000, i.e.~three orders of magnitude, regarding the online simulation time and of up to 30 regarding all necessary computing effort was observed in comparison to the conventional FE\textsuperscript{2} scheme. It was also empirically shown that the clustered training strategy not only reduces the training effort, but lowers the number of modes necessary in the online phase.

Comparing the obtained efficiency improvement measured by the speedup factor w.r.t.~the conventional FE\textsuperscript{2} method with those reported in the literature, at least four factors should be taken into account. Firstly, in some publications the speedup factor relates to RVE simulations, not actual multiscale simulations. Secondly, the complexity of the material formulations differs considerably among the papers. Thirdly, the publications use different accuracy requirements and fourthly, the redundancy of the RVEs' geometries and mesh densities must be assessed. The last factor is the most important one, because the main effect of the ROM method is to \enquote{filter out} redundancies, of which more are expected, if the geometry repeats in some way and the mesh density is very high. In the works of \citet{Raschi2021}, for example, the speedup of a composite RVE with different fiber stackings fell in the range of about $10$ to $10^4$, depending on the redundancy in the geometry. Exemplary, \citet{Rocha2020} obtained a speedup factor of 270 and \citet{Caicedo2019} report a speedup factor of 170--302, depending on the accuracy requirement. Other papers will come to comparable results, but the actual values strongly depend on the aforementioned factors.

The combination of the presented methods in this paper lie in terms of efficiency in the order of those observed in literature. However in many publications the offline effort is not quantified properly and the simulations settings differ substantially and are therefore difficult to be compared. Further research will show, which methods should be combined in the most effective way, to get robust and reliable simulation results in the multiscale context, while assuring low numerical costs and preserving most of the flexibility and modularity of the conventional FE\textsuperscript{2} method. Especially NN methods draw a lot of attention, since substantial online simulation cost reductions can be gained, e.g.~\cite{Le_neural_network_2015} found a speedup factor of a data driven NN approach w.r.t.~conventional FE\textsuperscript{2} of $120$,~in \citep{Huetter2021} a speedup of $142$ for a hybrid NN approach is mentioned, \citet{Logarzo2021} investigated a factor of $235$ for a \enquote{smart} data driven approach, \citet{Liu2019} reported a speedup of $930$ for a deep material network and \citet{Ghavamian2019} a speedup of three orders of magnitude for a recurrent NN, while due to the high training effort it is argued that, for a 		
fictitious example, a \enquote{break even} of over 100 macro simulations is needed until the training costs are compensated, where it must further considered that the recurrent NN looses its high predictive quality, as soon as it has to extrapolate to unseen data. \citet{Lu2019} report a speedup in the order of four magnitudes for a data-driven NN approach in the context of electrical homogenization. We conclude that in light of the high offline costs of NN methods in general, and the adaptation effort of hybrid NN approaches in particular, the high intrinsic adaptivity, ensured physical plausibility and low offline costs of the clustered hyper ROM method let it appear as an advantageous alternative to NN or data-driven techniques for many practical applications.

\section*{Declaration of competing interest}
The authors declare that they have no known competing financial interests or personal relationships that could have appeared to influence the work reported in this paper.
\section*{Acknowledgment}
The authors gratefully acknowledge computing time on the Compute Cluster of the Faculty of Mathematics and Computer Science of Technische Universität Bergakademie Freiberg, operated by the University Computing Center (URZ) and funded by the Deutsche Forschungsgemeinschaft (DFG) under DFG grant number 397252409.\\ 

\indent We thank our colleagues George Chatzigeorgiou and Fodil Meraghni (Arts et M\'etiers Institute of Technology, CNRS, Universit\'e de Lorraine) for supplying the data for the woven composite example.\\

\indent We are further grateful to Joaquin Hern\'andez (Centre Internacional de M\`etodes Num\`erics en Enginyeria (CIMNE), Technical University of Catalonia) for the fruitful discussions about details of his ECM method.
\appendix
\section{Derivation of the Algorithmically Consistent Tangent}
To get the algorithmic consistent tangent with respect to the Cauchy stress and the rate of deformation as needed by FE codes using the Updated Lagrangian formulation, approximations must be made, because at the material point level usually the relation $\tens[n+1]{\dot{\MStretch}}=\tens[n+1]{\dot{\MStretch}}(\tens[n+1]{\MRateofdef})$ is not known since $\tens[n+1]{\MRateofdef}$ is computed using an approximation (usually midpoint configuration after Hughes and Winget) which is only known at the element level. The quality of the approximation must be confirmed empirically. The tangent can be found by applying the product rule
\begin{equation}\label{tangent_continuum}
	\dfracderiv{\tens{\MStress}}{\tens{\MRateofdef}}=\dfracderiv{\tens{\MStress}}{\tens[\Biot]{\MStress}}:\dfracderiv{\tens[\Biot]{\MStress}}{\tens{\dot{\MStretch}}}:\dfracderiv{\tens{\dot{\MStretch}}}{\tens{\MRateofdef}}=\dfracderiv{\tens{\MStress}}{\tens[\Biot]{\MStress}}:\undertilde{\undertilde{C}}^\Biot:\dfracderiv{\tens{\dot{\MStretch}}}{\tens{\MRateofdef}}, \ \ \ \text{with} \ \ \ \undertilde{\undertilde{C}}^\Biot=\dfracderiv{\tens[\Biot]{\MStress}}{\tens{\dot{\MStretch}}}\ .
\end{equation}
\noindent Thereby $\dfracderiv{\tens{\MStress}}{\tens[\Biot]{\MStress}}$ is computed from equation \eqref{biotstress} as
\begin{equation}\label{d_Cauchy_d_Biot}
	\dfracderiv{\tens{\MStress}}{\tens[\Biot]{\MStress}}=\dfrac{2}{\MJacobian}\invers{[\indexnot{\MRotation}{ik}\invers{\indexnot{\MDefgrad}{lj}}+\invers{\indexnot{\MDefgrad}{ki}}\indexnot{\MRotation}{jl}]}\, \underline{E}_i\otimes\underline{E}_j\otimes\underline{E}_k\otimes\underline{E}_l
\end{equation}
\noindent From the RVE's FE system matrix $\undertilde{\undertilde{C}}^\Biot$ can be condensed out as shown in equation \eqref{Biot_Tangent}. The last term of equation \eqref{tangent_continuum} can only be approximated, whereby the following approximation, which assumes a constant rotation over the time step ($\tens{\dot{\MRotation}}\approx\tens{0}$) and thereby corresponding to the idea of the Updated Lagrangian Scheme, yields good results
\begin{equation}\label{dUdot_dd}
   \tens{\dot{\MStretch}}=\transp{\tens{\MRotation}}\cdot[\tens{\MRateofdef}+\tens{\MSpin}]\cdot\tens{\MDefgrad}-\transp{\tens{\MRotation}}\cdot\tens{\dot{\MRotation}}\cdot\tens{\MStretch}\ ,\ \ \ \dfracderiv{\tens{\dot{U}}}{\tens{\MRateofdef}}\approx\indexnot{\MRotation}{ki}\indexnot{\MDefgrad}{lj}\, \underline{E}_i\otimes\underline{E}_j\otimes\underline{E}_k\otimes\underline{E}_l.
\end{equation}
\noindent Equations \eqref{tangent_continuum}--\eqref{dUdot_dd} together yield
\begin{equation}
	\dfracderiv{\tens{\MStress}}{\tens{\MRateofdef}}\approx \dfrac{2}{\MJacobian}\invers{[\indexnot{\MRotation}{im}\invers{\indexnot{\MDefgrad}{nj}}+\invers{\indexnot{\MDefgrad}{mi}}\indexnot{\MRotation}{jn}]}C^{\Biot}_{mnop}\,\indexnot{\MRotation}{ko}\indexnot{\MDefgrad}{lp}\, \underline{E}_i\otimes\underline{E}_j\otimes\underline{E}_k\otimes\underline{E}_l\ .
\end{equation}
%In the actual implementation computing the appearing inverses should be avoided. At first the (symmetric) eigenvalue decomposition of the right Cauchy green tensor $\transp{\tens{\MDefgrad}}\cdot\tens{\MDefgrad}$ has to be computed. Having this decomposition, $\tens{\MStretch}$ and $\invers{\tens{\MStretch}}$ can easily be calculated and therewith $\tens{\MRotation}$ by a right multiplication of $\tens{\MDefgrad}$ with $\invers{\tens{\MStretch}}$. Having the inverse of the right stretch and the Rotation, the inverse of the deformation gradient can be computed by $\invers{\tens{\MDefgrad}}=\invers{\tens{\MStretch}}\cdot\transp{\tens{\MRotation}}$ . The inverse of equation \eqref{d_Cauchy_d_Biot} should never be computed actually, instead its matrix (in voigt notation) should be factorized to a LU decomposition and then the system of equations should be solved for the right hand side in terms of the stress resp.~tangent.
\label{sec:appendix_biot_tangent}
\bibliographystyle{unsrtnat}
\bibliography{literatur}
\end{document}